\documentclass[a4paper,12pt,nosumlimits]{amsart}
\usepackage{latexsym,amsmath,amssymb}
\usepackage[english]{babel}
\usepackage{amscd}
\title[Gauging away noncommutativity]
{Moduli spaces of noncommutative instantons: \\[6pt]
gauging away noncommutative parameters}
\date{v1: 24 September 2009; v2: 30 September 2010} 
\author{Simon Brain}
\address{\flushleft SISSA, Via Beirut 2-4, 34151 Trieste, Italia}
\email{brain@sissa.it}
\author{Giovanni Landi}
\address{\flushleft Dipartimento di Matematica e
Informatica, Universit\`{a} di Trieste, Via A. Valerio 12/1, 34127
Trieste, Italia, and INFN, Sezione di Trieste, Trieste, Italia}
\email{landi@univ.trieste.it}
 \numberwithin{equation}{section}

\newcommand{\half}{ \tfrac{1}{2} }

\newtheorem{thm}{Theorem}[section]

\newtheorem{lem}[thm]{Lemma}
\newtheorem{prop}[thm]{Proposition}
\newtheorem{rem}[thm]{Remark}

\theoremstyle{definition}
\newtheorem{defn}[thm]{Definition}
\newtheorem{example}[thm]{Example}

\newcommand{\pp}{{\sf p}}
\newcommand{\qp}{{\sf q}}
\newcommand{\sfu}{{\sf u}}

\newcommand{\sfm}{\sf m}
\newcommand{\sfM}{\sf M}
\newcommand{\sfV}{{\sf V}}
\newcommand{\Pp}{{\sf P}}
\newcommand{\Qp}{{\sf Q}}

\newcommand{\ii}{\mathrm{i}}
\newcommand{\E}{\mathcal{E}}

\newcommand{\A}{\mathcal{A}}
\newcommand{\mH}{\mathcal{H}}
\newcommand{\mK}{\mathcal{K}}
\newcommand{\mL}{\mathcal{L}}
\newcommand{\mM}{\mathcal{M}}
\newcommand{\tmM}{\widetilde{\mathsf{M}}}
\newcommand{\mC}{\mathcal{C}}
\newcommand{\mR}{\mathcal{R}}
\newcommand{\B}{\mathcal{B}}

\newcommand{\hatA}{\widehat{A}}
\newcommand{\la}{\langle}
\newcommand{\ra}{\rangle}

\newcommand{\tr}{\triangleright}
\newcommand{\tl}{\triangleleft}

\newcommand{\n}{\nabla}
\newcommand{\ep}{\epsilon}

\newcommand{\M}{\textup{M}}
\newcommand{\SL}{\textup{SL}}
\newcommand{\Sp}{\textup{Sp}}

\newcommand{\SU}{\textup{SU}}

\newcommand{\HH}{\mathbb{H}}
\newcommand{\ZZ}{\mathbb{Z}}
\newcommand{\C}{\mathbb{C}}
\newcommand{\TT}{\mathbb{T}}

\newcommand{\RR}{\mathbb{R}}

\newcommand{\End}{\textup{End}}

\newcommand{\id}{\textup{id}}
\newcommand{\D}{\textup{d}}
\newcommand{\ch}{\textup{ch}}
\newcommand{\utimes}{\,\underline{\otimes}\,}

\renewcommand{\o}{{}_{\scriptscriptstyle(1)}}
\renewcommand{\t}{{}_{\scriptscriptstyle(2)}}
\newcommand{\thr}{{}_{\scriptscriptstyle(3)}}
\newcommand{\fo}{{}_{\scriptscriptstyle(4)}}

\newcommand{\bo}{{}^{{\scriptscriptstyle(-1)}}}
\newcommand{\bt}{{}^{{\scriptscriptstyle(0)}}}
\newcommand{\bp}{{}^{{\scriptscriptstyle(+1)}}}

\newcommand{\uo}{{}_{\scriptscriptstyle(\underline 1)}}
\newcommand{\ut}{{}_{\scriptscriptstyle(\underline 2)}}

 \def\lcross{{>\!\!\!\tl}}

\def\lcocross{{>\!\!\blacktriangleleft}}
\def\rbiprod{{\cdot\kern-.31em\tr\!\!\!<}}
\def\lbiprod{{>\!\!\!\triangleleft\kern-.33em\cdot}}
\def\dcross{{\bowtie}}

\begin{document}
\begin{abstract}
Using the theory of noncommutative geometry in a braided monoidal
category, we improve upon a previous construction of noncommutative families of
instantons of arbitrary charge on the deformed sphere $S^4_\theta$.
We formulate a notion of noncommutative parameter spaces for
families of instantons and we explore what it means for such families
to be gauge equivalent, as well as showing how to remove gauge parameters using
a noncommutative quotient construction. Although the parameter spaces are
{\em a priori} noncommutative, we show that one may always recover a
classical parameter space by making an appropriate choice of gauge
transformation.\end{abstract}

\maketitle

\tableofcontents
\section{Introduction}
\label{section introduction} We study families of instantons on the
noncommutative four-sphere $S^4_\theta$ of \cite{cl:id} and what it
means for such families to be gauge equivalent. We show that,
although it is perfectly natural to allow for the possibility of
families of instantons parameterised by noncommutative spaces, these
`noncommutative parameters' may always be removed by an appropriate
choice of gauge transformation so as to recover a `commutative'
parameter space.

The study of instantons on $S^4_\theta$ was initiated in
\cite{lvs:pfns,lvs:nitcs} and developed further in \cite{lprs:ncfi,
bl:adhm}, where it was observed that one may construct families of
instantons which are parameterised by noncommutative spaces. These
parameter spaces arise in a very natural way and suggest that we
should consider seriously the idea that the moduli space of
instantons might also be noncommutative. On the other hand,
instantons on $S^4_\theta$ are defined in terms of absolute minima
of the Yang-Mills energy functional; thinking of the moduli space as
being `the set of all such minima modulo gauge equivalence'
naturally leads us to expect it to be a classical space. Our goal in
the present article is to use gauge theory to reconcile this
apparent dichotomy between classical and noncommutative parameter
spaces.

The driving force behind our investigation is the fact that the
quantum sphere $S^4_\theta$ can be obtained from its classical
counterpart $S^4$ by means of a Hopf cocycle twisting procedure. The
often-forgotten feature of this construction is that it deforms not
just the four-sphere but in fact the entire category in which it
lives. In our case, with $H=\A(\TT^2)$ the Hopf algebra of
coordinate functions on the two-torus $\TT^2$, the deformation takes
the form of a `quantisation functor' from the category of
$H$-comodules, wherein lives the classical sphere $S^4$, to a new
category containing the quantum sphere $S^4_\theta$. This new
category is the category of comodules for a twisted Hopf algebra
$H_F$, with $F$ a twisting Hopf cocycle.

By expressing the construction of instantons on $S^4$ entirely in
this categorical framework, we are able to apply the quantisation
functor and hence obtain a construction of instantons on
$S^4_\theta$. Since all parameter spaces we consider are themselves
objects in the category, they are twisted as well by the functor and
we are naturally led to the concept of noncommutative families of
instantons. We discuss what it means for such families to be gauge
equivalent and show, just as in the classical case, how one can
quotient parameter spaces by the resulting equivalence relation. As
mentioned, a suitable choice of gauge transformation can be used to
remove the noncommutativity of the parameters and produce an
equivalent description in terms of usual spaces.

The paper is organised as follows. After section
\S\ref{se:hap}, which reviews the abstract theory of Hopf algebras
and the cocycle twisting construction, we give an overview of how to
construct the various noncommutative spaces that we shall need. In
particular, \S\ref{section noncomm hopf} recalls the construction of
the noncommutative $\SU(2)$ Hopf fibration $S^7_\theta\to
S^4_\theta$ using cocycle twisting, together
with the canonical differential structures on these spaces.

As a way to understand the structures involved, the first topic of
the paper will be to study parameter spaces for charge one
instantons. In the classical case, one can construct all such
instantons by acting upon a basic instanton with the group
$\SL(2,\HH)$ of conformal transformations of the four-sphere $S^4$.
In \S\ref{section quantum conformal group} and \S\ref{sect braided}
we write the various symmetry groups of $S^4$ in an entirely
$H$-covariant setting, which we then twist using the quantisation
functor. This leads naturally to `braided geometry' in the deformed
category and, in particular, to a braided Hopf algebra
$\B(\SL_\theta(2,\HH))$ of conformal symmetries. It obeys the usual
axioms of a Hopf algebra, but with its structure maps required to be
morphisms in the category. To this we apply a cobosonisation process
to recover an `ordinary' Hopf algebra, which takes the form of a
Hopf algebra biproduct $\B(\SL_\theta(2,\HH))\lbiprod H_F$.

In \S\ref{section gauge theory} we review the basic notions of gauge
theory on $S^4_\theta$, then generalise them by formulating a notion
of noncommutative parameter spaces for families of instantons and
what it means for such families to be gauge equivalent. Using these
definitions, we are able to parallel the classical case by
interpreting the quantum symmetry group
$\B(\SL_\theta(2,\HH))\lbiprod H_F$ as a parameter space for the set
of charge one instantons on the quantum four-sphere. In \S\ref{sect
charge one params} we study this noncommutative parameter space in
more detail, seeking where possible to remove all parameters
corresponding to gauge equivalence instantons. In the classical
case, the gauge parameters are described by the subgroup $\Sp(2)$ of
$\SL(2,\HH)$ consisting of isometries of the sphere, so the `true'
parameter space is the quotient $\SL(2,\HH)/\Sp(2)$. In the
noncommutative case there is a braided group of isometries
$\B(\Sp_\theta(2))$ that we are immediately able to remove by means
of a quantum quotient construction.

Far more subtle is the question of how to remove the gauge
parameters corresponding to the subalgebra $H_F$ of
$\B(\SL_\theta(2,\HH))\lbiprod H_F$. These extra symmetries
correspond to the inner automorphisms of the coordinate algebra
$\A(S^4_\theta)$ of the deformed 4-sphere and constitute a very important
part of its noncommutative geometry \cite{ac:fncg}. We
show that there are many ways in which to quotient $H_F$ away from
$\B(\SL_\theta(2,\HH))\lbiprod H_F$; in the classical case, every
way we do this gives the same answer, but in the noncommutative case
we get families of parameter spaces which are clearly different
(some being quantum, some classical) but all have the same classical
limit. We show that these parameter spaces are all gauge equivalent,
finding as special cases both a commutative parameter space as well
as the noncommutative parameter space found previously in
\cite{lprs:ncfi}.

In \S\ref{sect adhm} we see how this method generalises to
instantons with higher charge. We review the usual ADHM construction
of \cite{adhm:ci} ({\em cf.} also \cite{ma:gymf}) in the context of
braided geometry, which we then deform using the quantisation
functor. With a few minor differences this essentially reproduces
the noncommutative ADHM construction of \cite{bl:adhm}, although
derived from a different and arguably more natural approach. As in
the charge one case, we show how to remove the gauge parameters
corresponding to the torus algebra $H_F$, finding in particular that
a certain choice of gauge yields again a commutative parameter
space.

The paper concludes with an appendix reviewing the
notion of quantum families of maps, which is an essential theme used
throughout the paper. In looking for moduli spaces of instantons,
our philosophy is to look not for a set of objects but rather for a
space which parameterises those objects, that is to say we ask for
some geometric structure. In categorical terms, this means defining
a functor from the category of algebras to the category of sets and
then looking for the moduli space as a universal object; this is
necessarily an object in the source category, \emph{i.e.} an
algebra.

\section{Preliminaries on Hopf Algebras and their Deformations}\label{se:hap}

We review here some important elements of Hopf algebra theory,
including the cocycle twisting construction that will play such an
important part in what follows.

\subsection{Hopf algebra preliminaries} \label{section hopf algebra prelims}
We recall some basic facts from the theory of Hopf algebras and
related structures following mainly \cite{ma:book} ({\em cf}. also
\cite{ma:habc}). Given a Hopf algebra $H$ over $\C$ we denote its
coproduct, counit and antipode by $\Delta:H \rightarrow H\otimes H$,
$\ep:H\rightarrow \C$ and $S:H\rightarrow H$, respectively. The product map is
usually suppressed, although when explicitly written it is denoted
$m(g \otimes h)=gh$. We use Sweedler notation for the coproduct,
$\Delta h=h\o \otimes h\t$; also we indicate $(\Delta \otimes \id) \circ \Delta h = (\id \otimes \Delta) \circ \Delta h = h\o \otimes h\t \otimes h \thr$ and so on, with summation inferred.  A Hopf algebra
$H$ is said to be {\em coquasitriangular} if it is equipped with a
convolution-invertible Hopf bicharacter $\mathcal{R}:H\otimes
H\rightarrow \C$ satisfying
$$g\o h\o\mathcal{R}(h\t,g\t)=\mathcal{R}(h\o,g\o)h\t g\t$$ for all
$g,h\in H$. Convolution invertibility is the existence of a map
$\mathcal{R}^{-1}:H\otimes H\rightarrow \C$ such that
\begin{equation}\label{con-inv}
\mathcal{R}(h\o,g\o)\mathcal{R}^{-1}(h\t,g\t)=\mathcal{R}^{-1}(h\o,g\o)\mathcal{R}(h\t,g\t)=\ep(g)\ep(h)
\end{equation}
for all $g,h \in H$. On the other hand being a bicharacter means that
\begin{equation}\label{hopf-bic}
\mathcal{R}(fg,h)=\mathcal{R}(f,h\o)\mathcal{R}(g,h\t), \qquad \mathcal{R}(f,gh)=\mathcal{R}(f\o,h)\mathcal{R}(f\t,g)
\end{equation}
for all $f,g,h \in H$. If in addition $\mathcal{R}$ obeys the
identity $$\mathcal{R}(b\o,a\o)\mathcal{R}(a\t,b\t)=\ep(a)\ep(b)$$
for all $a,b\in H$ then we say that $H$ is a {\em cotriangular} Hopf
algebra.

A left module structure $H\otimes A \to A$ on a vector space $A$ is
denoted $\tr$, {\em i.e.} we write $h \otimes a \mapsto h\tr a$ for
$h \in H$, $a \in A$. A right module structure is denoted $\tl$.
Similarly we denote a left comodule structure on $A$ by $\Delta_L:A
\rightarrow H \otimes A$, again using Sweedler notation
$\Delta_L(a)=a\bo \otimes a\bt$. A right comodule structure is
written $\Delta_R:A\rightarrow A \otimes H$ with a similar Sweedler
notation: $\Delta_R(a)=a\bt \otimes a\bp$. We denote the categories
of left $H$-modules and left $H$-comodules by ${}_H\mM$ and
${}^H\mM$ respectively. Moreover, we say $A$ is a {\em left crossed
$H$-module} if it is both a left $H$-module and a left $H$-comodule
and these structures obey the compatibility condition
\begin{equation}\label{eqn cross module condition} h\o a\bo\otimes h\t\tr a\bt=(h\o\tr a)\bo h\t\otimes (h\o \tr a)\bt\end{equation} for all $h \in H$ and $a \in A$. The category of
left crossed $H$-modules is denoted ${}_H^H\mC$, with a similar
definition for the category $\mC_H^H$ of right crossed $H$-modules.
When $H$ is coquasitriangular, there is a canonical monoidal functor
${}^H\mM\rightarrow {}^H_H\mC$ given by equipping an $H$-comodule
$A$ with the $H$-action
\begin{equation}\label{can act}
h\tr a:=\mathcal{R}(a\bo,h)a\bt, \qquad a \in A,~ h \in H,
\end{equation}
where $a \mapsto a\bo\otimes a\bt$ denotes the $H$-coaction, as before. One may
check that this gives a well-defined $H$-action using the fact that
$\mathcal{R}$ is a Hopf bicharacter.

A monoidal (or tensor) category 
is {\em braided} if for each pair of objects $V,W$ there is an
isomorphism $\Psi_{V,W}:V\otimes W\rightarrow W\otimes V$, obeying
certain natural hexagon identities \cite{js}. The simplest example is the
category $\mathsf{Vec}$ of complex vector spaces, with the
monoidal structure given by the usual tensor product of vector
spaces and braided by the flip map:
$\Phi_{V,W}(v\otimes w)=w\otimes v$ for all $v \in V$ and $w \in W$. More
generally, if $H$ is a Hopf algebra then the category ${}^H\mM$ has
a monoidal structure given by the tensor product coaction,
\begin{equation}\label{tpco}
\Delta_{V\otimes W}(v\otimes w)=v\bo w\bo \otimes v\bt\otimes
w\bt, \qquad v\in V, ~w \in W.
\end{equation}
If in addition $H$ is
coquasitriangular then ${}^H\mM$ is braided by the collection of
morphisms
\begin{equation}\label{bramor}
\Psi_{V,W}(v\otimes w)=\mathcal{R}(w\bo,v\bo)w\bt\otimes v\bt ,
\end{equation}
for each pair $V$, $W$ of
left $H$-comodules with $v\in V$ and $w \in W$. In particular, if $A$
and $B$ are left $H$-comodule algebras ({\em i.e.} algebras in the
category ${}^H\mM$),  the braiding allows one to give
a tensor product algebra $A\underline{\otimes}B$,
with the product
$$
(a\otimes b)(c\otimes d)=a \Psi_{B,A}(b\otimes c)d ,
$$
and which lives in the category ${}^H\mM$ by the
coaction $\Delta_{A\otimes B}$ in \eqref{tpco} above.  The symbol
$\underline{\otimes}$ is to remind us that the tensor product is the
braided one. The braided monoidal category
$(\mathsf{Vec},\otimes,\Phi)$ is recovered by putting $H=\C[\C]$,
the coordinate algebra of the complex numbers, with its trivial
coquasitriangular structure.

If $A$ is a left $H$-module algebra ({\em i.e.} an algebra in the
category ${}_H\mM$), then there is a cross product
algebra $A \lcross H$ built on $A \otimes H$ as a vector space, with
algebra structure
\begin{equation}\label{eqn cross product}
(a \otimes g)(b \otimes h):=a(g\o\tr b) \otimes g\t h
\end{equation}
for all $g,h \in H$ and $a,b \in A$; and unit $1_H\otimes 1_A$. Equally well, if $A$ is a
left $A$-comodule coalgebra ({\em i.e.} a coalgebra in the category
${}^H\mM$), there is a cross coproduct coalgebra $A\lcocross H$
built on $A \otimes H$ as a vector space, with counit $\ep=\ep_{_A} \otimes \ep_{_H}$ and coproduct
\begin{equation}\label{eqn cross coproduct}
\Delta(a \otimes h):= a\o\otimes a\t\bo h\o\otimes a\t\bt \otimes h\t\end{equation}
for all $h \in H$ and $a \in A$.

Furthermore, we may consider bialgebras and Hopf algebras which are
themselves objects in a braided category. A bialgebra in the
category ${}^H\mM$ is by definition a bialgebra in the usual sense,
{\em i.e.} it obeys all of the usual axioms, but with its structure
maps now as morphisms in the category. We call such an object $K$ a
{\em braided bialgebra}; we denote its structure maps by
($\underline{m}$, $\underline{\Delta}$, $\underline{\ep}$) if we wish
to stress that they are morphisms in a braided category. In
particular  the coproduct $\underline{\Delta}:K\to
K\underline{\otimes}K$ is required to be an algebra homomorphism
from $K$ into the braided tensor product. If $K$ has also an
antipode $\underline{S}$ (obeying the usual axioms, but again
required to intertwine the $H$-coaction) then we say that $K$ is a
{\em braided Hopf algebra}.

If $K$ is a braided bialgebra in ${}^H\mM$, we already mentioned that it becomes an
$H$-module algebra {\em via} the canonical action \eqref{can act}. It follows that
the vector space $K\otimes H$ may be equipped with the structure of
an `ordinary' bialgebra, given by the above cross product and cross
coproduct constructions in \eqref{eqn cross product} and \eqref{eqn cross coproduct} respectively.
The resulting bialgebra is called the {\em cobosonisation} of $K$ and is denoted $K \lbiprod H$. A sufficient condition for $K \lbiprod H$ to be a Hopf algebra is that $K$ and $H$
be Hopf algebras with the antipode of $H$ invertible, in which case
the antipode of $K \lbiprod H$ is
\begin{equation}\label{eqn biproduct antipode}
S(a\otimes h)=(1\otimes S^{-1}(a\bo h))(\underline{S}(a\bt)\otimes 1).
\end{equation}
The cobosonisation is of special interest because left $K \lbiprod
H$-comodules are by construction exactly the same thing as left
$K$-comodules in the category ${}^H\mM$. The cobosonisation thus
gives us the option of working either with a braided Hopf algebra
$K$ and its braided comodules or simply with comodules for the
ordinary Hopf algebra $K \lbiprod H$.

\subsection{Cocycle deformations of Hopf algebras}\label{section cocycle twists}
We turn now to recalling the basic theory of
quantisation of Hopf algebras by `cocycle cotwist' as in \cite{ma:book}.
Starting with a two-cocycle $F$ on a Hopf algebra $H$, we give the
relevant formul{\ae} for obtaining a deformed Hopf algebra $H_F$ as
well as the appropriate deformations of the algebras or coalgebras
on which $H$ acts or coacts. We illustrate the theory with the
well-known example of the noncommutative $n$-torus that we shall
also use later on in the paper.

By a {\em two-cocycle} on a Hopf algebra $H$ we mean a map $F:H
\otimes H\rightarrow \C$ which is unital, convolution-invertible in
the sense of \eqref{con-inv} and obeys the cocycle condition $\partial F=1$ or
\begin{equation} \label{eqn two cocycle}
F(g\o,f\o)F(h\o,g\o f\t)F^{-1}(h\t g\thr,f\thr)F^{-1}(h\thr,g\fo)=\ep(f)\ep(h)\ep(g)
\end{equation}
for all $f,g,h \in H$.
Given such an $F$, there is a {\em cotwisted}
Hopf algebra $H_F$ which as a coalgebra is the same as $H$ but whose
product is replaced by
\begin{equation} \label{eqn twisted product}h \bullet_{_F}
g=F(h\o,g\o)h\t g\t F^{-1}(h\thr,g\thr)\end{equation} and whose
antipode becomes
\begin{equation} \label{eqn twisted antipode}
S_F(h):=U(h\o)S(h\t)U^{-1}(h\thr), \qquad \textup{with} \qquad U(h):=F(h\o,Sh\t).
\end{equation}
The cocycle condition \eqref{eqn two cocycle} assures that the product in $H_F$ is associative.
If $H$ has a coquasitriangular
structure $\mathcal{R}:H\otimes H\rightarrow \C$, then $H_F$ is also
coquasitriangular with
\begin{equation} \label{eqn twisted
R-matrix}\mathcal{R}_F(h,g):=F(g\o,h\o)\mathcal{R}(h\t,g\t)F^{-1}(h\thr,g\thr).
\end{equation}
In the case where $H$ is commutative, then $\mathcal{R}_F$ in fact
defines a cotriangular structure on $H_F$ and, as a consequence, the
induced braiding $\Psi$ on the category ${}^{H_F}\mM$ is symmetric
in the sense that $\Psi^2=\textup{id}$.

In passing from $H$ to $H_F$ one finds that ${}^H\mM$ and
${}^{H_F}\mM$ are isomorphic as braided monoidal categories. Indeed,
since the cotwist does not change the coalgebra structure of $H$, it
follows that $H$-comodules are also $H_F$-comodules and $H$-comodule
morphisms are also $H_F$-comodule morphisms; thus there is a functor
$\mathcal{G}_F:{}^H\mM \to {}^{H_F}\mM$ which leaves the coactions
unchanged. As categories, we have simply that ${}^H\mM$ and
${}^{H_F}\mM$ are just the same: the non-trivial part of the
isomorphism is contained in what happens to the monoidal structure.
Writing $V_F:=\mathcal{G}_F(V)$, the category ${}^{H_F}\mM$ gets a
new monoidal structure by
$$\sigma_F:V_F\otimes W_F\to
(V\otimes W)_F, \qquad v\otimes w \mapsto F(v\bo,w\bo)v\bt\otimes
w\bt.$$ One checks \cite{mo:twist} that $\mathcal{G}_F$ is a
monoidal functor and that it intertwines the braidings in ${}^H\mM$
and ${}^{H_F}\mM$ given respectively by $\mathcal{R}$ and
$\mathcal{R}_F$ according to the formula \eqref{bramor}. We call
$\mathcal{G}_F$ the `quantisation functor' associated to the cocycle
$F$ since it simultaneously deforms all $H$-covariant constructions
to corresponding versions which are covariant under $H_F$.

In particular, if $A$ is an algebra in the category ${}^H\mM$, then
under the functor $\mathcal{G}_F$ the product map $m:A\otimes A\to
A$ becomes a map $(A\otimes A)_F \to A_F$. Composing this with
$\sigma_F$ yields a new product map
\begin{equation}\label{eqn twisted prod}m_F:A_F\otimes
A_F\to A_F, \qquad a\otimes b \mapsto a\cdot_F b:=F(a\bo,b\bo)a\bt
b\bt,\end{equation} and $m_F$ automatically makes $A_F$ into an
$H_F$-comodule algebra.

In the same way, if $A$ is a coalgebra in the category ${}^H\mM$ with coproduct
$\Delta:A\to A\otimes A$,  applying the functor
$\mathcal{G}_F$ results in a map $A_F\to (A\otimes A)_F$. Then,
composing with $\sigma_F^{-1}$ yields a new coproduct map
\begin{equation}\label{eqn twisted coprod}\Delta_F:A_F\to
A_F\otimes A_F, \qquad a \mapsto F^{-1}(a\o{}\bo,a\t{}\bo)a\o{}\bt
\otimes a\t{}\bt,\end{equation} which automatically makes $A_F$ into
an $H_F$-comodule coalgebra with counit $\underline{\ep}_F=
\underline{\ep}$.

We may of course put these two constructions together. Suppose that
$A$ is a bialgebra in the category ${}^H\mM$, {\em i.e.} it is both
an $H$-comodule algebra and an $H$-comodule coalgebra in a
compatible way. Then we may simultaneously twist the product and
coproduct on $A$ and, as one might expect \cite{ma:qst}, the result is a bialgebra
$A_F$ in the category ${}^{H_F}\mM$. Moreover, if $A$
is a Hopf algebra in ${}^H\mM$ with antipode $\underline{S}$, then
$A_F$ is a Hopf algebra in ${}^{H_F}\mM$ with antipode
$\underline{S}_F=\underline{S}$.

\begin{rem}\label{rem real cocycle}
\textup{
In the case where $H$ is a Hopf $*$-algebra, (thus in particular
$\Delta$ is a $*$-algebra map with $(S \circ *)^2=\textup{id}$), we need to add the condition that $F$ is a {\em real cocycle} in
the sense that
\begin{equation} \label{eqn real cocycle} \overline{F(h,g)}=F((S^2g)^*,(S^2h)^*).\end{equation} Then $H_F$
acquires a deformed $*$-structure
\begin{equation} \label{eqn twisted star} h^{*_F}:=\overline{V^{-1}(S^{-1}h\o)}(h\t)^*\overline{V(S^{-1}h\thr)},
\qquad \textup{with} \qquad
V(h):=U^{-1}(h\o)U(S^{-1}h\t).\end{equation} Then, if $A$ is a left
$H$-comodule algebra and a $*$-algebra such that the coaction is a
$*$-algebra map, the twisted algebra $A_F$ gets a new $*$-structure
as well,
\begin{equation} \label{eqn twisted comodule star}
a^{*_F}:=\overline{V^{-1}(S^{-1}a\bo)}(a\bt)^*.\end{equation}
}
\end{rem}

\begin{example}\label{example twisted torus}
The Hopf algebra $H:=\A(\TT^n)$ of functions on the $n$-torus $\TT^n$
is the algebra
\begin{equation}H:=\A[t_j,t_j^{-1} ~|~j=1,\ldots,n]\end{equation} equipped with the
Hopf $*$-algebra structure
\begin{equation}t_j^*=t_j^{-1}, \quad \Delta(t_j)=t_j \otimes t_j,
\quad \ep(t_j)=1, \quad S(t_j)=t_j^{-1}\end{equation} for all
$j=1,\ldots, n$, and, as usual, $\Delta$, $\ep$ extended as
$*$-algebra maps and $S$ extended as an anti-$*$-algebra map. The
canonical right action of $\TT^n$ on itself by group multiplication
dualises to give a left coaction
\begin{equation}\label{grouplike}\Delta_L:\A(\TT^n)\rightarrow
H\otimes\A(\TT^n), \qquad u_j \mapsto t_j\otimes u_j,\end{equation}
where we write $u_j$, $u_j^*$, $j=1,\ldots,n$, for the generators of
$\A(\TT^n)$ viewed as a left comodule algebra over itself. This
coaction is equivalent to a grading of $\A(\TT^n)$ by the Pontrjagin
dual group $\ZZ^n$ of $\TT^n$, for which the homogeneous elements
are the monomials of the form \begin{equation}\label{hom
grad}t^{\vec a}:=t_1^{a_1}t_2^{a_2}\cdots t_n^{a_n},\qquad \vec
a:=(a_1,a_2,\ldots,a_n)\in \ZZ^n.\end{equation} One defines a
two-cocycle $F$ on $H$ by choosing a real antisymmetric $n\times n$
matrix $\Theta=(\Theta_{jl})$ and setting
\begin{equation}\label{eqn torus cocycle}
F(t^{\vec a},t^{\vec b}):= \exp\left( \ii \pi (\vec a \cdot \Theta \cdot \vec b)\right)
\end{equation}
for homogeneous multi-degree elements $t^{\vec a},t^{\vec b}\in H$
and extended by linearity. It is straightforward to verify that $F$
is a cocycle which is real in the sense of Remark~\ref{rem real
cocycle}. Moreover, from its form as an exponential, this $F$ is a
Hopf bicharacter ({\em cf}. Eq.~\eqref{hopf-bic}), that is
$F(fg,h)=F(f,h\o)F(g,h\t)$ and $F(f,gh)=F(f\o,h)F(f\t,g)$ for all
$f,g,h \in H$. As a consequence it obeys $$F(Sh,g)=F^{-1}(h,g),
\quad F(h,Sg)=F^{-1}(h,g), \quad F(Sh,Sg)=F(h,g)$$ for all $g,h \in
H$. These properties mean that $F$ is determined by its values on
the generators $t_j$, for which we have
\begin{equation}\label{eqn cocycle on gens}
F(t_j,t_l)=\exp(\ii \pi \Theta_{jl}), \qquad j,l=1,\ldots,n.
\end{equation}
The product, $*$-structure and antipode on $H$ are in fact
undeformed by $F$, so $H=H_F$ as a Hopf $*$-algebra. However, the
trivial cotriangular structure $\mR=\ep \otimes \ep$ of $H$ twists
into
$$
\mR_F(t_j,t_l)=F(t_l,t_j)F^{-1}(t_j,t_l)=F^{-2}(t_j,t_l).
$$
As mentioned for the general construction, the deformation takes the
form of an isomorphism of braided monoidal categories from ${}^H\mM$
to ${}^{H_F}\mM$. In particular, for $\A(\TT^n)$, the effect is
that, considered as left $H$-comodule algebra for itself, the
$*$-structure on $\A(\TT^n)$ is unchanged but the product is twisted
into a new product:
$$
u_j \cdot_F u_l=u_j u_l F(t_j,t_l)=u_j u_l e^{ \ii\pi\Theta_{jl}}.
$$
We denote by $\A(\TT^n_\Theta)$ the $*$-algebra generated by the
$u_j$, $u_j^*$ with this new product; there are now relations
$$u_j \cdot_F u_l=e^{2\ii \pi \Theta_{jl}} u_l \cdot_F u_j, \qquad
u_l^*\cdot_F u_j=e^{2\ii\pi \Theta_{jl}} u_j\cdot_F u_l^*$$ for
each pair of indices $j,l=1,\ldots,n$. The original torus $\TT^n$
has been quantised to give the noncommutative torus $\TT^n_\Theta$.
\end{example}

\section{Hopf Fibrations over Spheres}\label{section noncomm hopf}
In this section we review the construction in \cite{lvs:pfns} of the
 $\SU(2)$-Hopf fibration over the noncommutative four-sphere $S^4_\theta$ of \cite{cl:id}.
 We begin by giving the
classical fibration in a coordinate algebra form, which we then
quantise by means of a cocycle cotwist. We then review how the same
twisting procedure also yields canonical differential calculi on the
noncommutative algebras, as well as a Hodge operator $*_\theta$ on
the sphere $S^4_\theta$.

\subsection{The classical Hopf bundle}The coordinate algebra $\A(\C^4)$ of the vector space
$\C^4$ is the associative commutative $*$-algebra generated by the
functions $z_j$, $j=1,\ldots,4$, together with their conjugates
$z_j^*$, $j=1,\ldots,4$. The coordinate algebra $\A(S^7)$ of the
seven-sphere is the quotient of $\A(\C^4)$ by the sphere relation
\begin{equation}\label{eqn seven sphere}
z_1^*z_1+z_2^*z_2+z_3^*z_3+z_4^*z_4=1.
\end{equation}

It is useful to arrange the generators of the algebra $\A(\C^4)$ into the matrix
\begin{equation}\label{u}
u:=\begin{pmatrix}z_1&z_2&z_3&z_4\\-z_2^*&z_1^*&-z^*_4&z_3^*\end{pmatrix}^{\textrm{t}},
\end{equation}
with ${}^{\textrm{t}}$ denoting matrix transposition, which we use
to give a right action of the classical group $\SU(2)$ on $\A(\C^4)$
by
$$
u \mapsto u w, \qquad \textup{with} \quad w=\begin{pmatrix}w^1 &
-\bar w^2\\ w^2& \bar w^1\end{pmatrix}\in \SU(2).
$$
This action preserves the sphere relation \eqref{eqn seven sphere}, whence it restricts to
an action of $\SU(2)$ on the coordinate algebra $\A(S^7)$ of the seven-sphere. The subalgebra
$\textrm{Inv}_{_{\SU(2)}}(\A(S^7))$ of invariant functions is
generated as a commutative $*$-algebra by the elements
\begin{equation}\label{eqn alg inc}\alpha=2(z_1z^*_3 + z_2^*z_4), \quad \beta=2(z_2z_3^* -
z_1^*z_4), \quad x= z_1z^*_1 + z_2z^*_2 - z_3z_3^* -
z_4z_4^*,\end{equation} together with their conjugates $\alpha^*$,
$\beta^*$ and $x^*=x$. It follows from the sphere relation
\eqref{eqn seven sphere} that these generators obey the relation
\begin{equation}\label{eqn four sphere}
\alpha^*\alpha+\beta^*\beta+x^2=(z_1^*z_1+z_2^*z_2+z_3^*z_3+z_4^*z_4)^2=1,
\end{equation} whence the
invariant subalgebra is the coordinate algebra of a four-sphere,
$$\textrm{Inv}_{\SU(2)}(\A(S^7))=\A(S^4).$$

Indeed, the sphere relation \eqref{eqn seven sphere} is equivalent
to requiring that $u^*u=1$, from which we automatically have that the
matrix-valued function
\begin{equation}\label{eqn basic proj} q:=uu^*=\tfrac{1}{2}\begin{pmatrix} 1+x & 0 & \alpha
& -\beta^* \\ 0 & 1+x & \beta & \alpha^* \\
\alpha^* & \beta^* & 1-x & 0 \\ -\beta & \alpha & 0 &
1-x\end{pmatrix}\end{equation} is a self-adjoint idempotent: $q^2=q=q^*$. Clearly
one has
$$(u w)(u w)^*=u (ww^*)u^* = uu^*, \qquad w \in \SU(2),
$$
and the entries of $q$ really do generate an $\SU(2)$-invariant
subalgebra. Moreover, Eq.~\eqref{eqn alg inc} defines an inclusion
of algebras $\A(S^4)\hookrightarrow\A(S^7)$, which is just a
coordinate algebra description of the standard Hopf fibration
$S^7\to S^4$ having  $\SU(2)$ as structure group.

\subsection{The noncommutative Hopf fibration}\label{section nc hopf fib}
We now obtain a noncommutative version of the Hopf fibration using
the method of `cocycle cotwisting' as described in \S\ref{section
cocycle twists}, with compatible torus (co)actions on the total and
the base spheres. The `deforming' Hopf algebra will be the algebra
$H:=\A(\TT^2)=\A[t_j,t_j^{*} ~|~j=1,2]$ of functions on the
two-torus $\TT^2$ with a `deforming' two-cocycle $F$ as given in
\eqref{eqn torus cocycle}. For the present case $\Theta$ is a real
$2\times 2$ antisymmetric matrix and hence of the form
$$
\Theta= \half \begin{pmatrix}0&\theta\\-\theta&0\end{pmatrix},
$$
with $\theta\in \RR$. We know from Example~\ref{example twisted
torus} that the twisted Hopf $*$-algebra structure on $H=\A(\TT^2)$
is in fact unchanged, so that $H=H_F$ as a Hopf $*$-algebra,
although the trivial cotriangular structure $\mR=\ep \otimes \ep$ of
$H$ is twisted into $\mR_F(t_j,t_l):=F^{-2}(t_j,t_l)$, leading to a
twisted product on $H$-comodule algebras.  Indeed, we now use $H_F$
to deform the Hopf fibration described in the previous section.

There is a left coaction of $H=\A(\TT^2)$ on the coordinate algebra
$\A(S^7)$ given by
\begin{eqnarray}\label{eqn coact spheres}
\Delta_L:\A(S^7)\to H\otimes\A(S^7), \qquad
\Delta_L(z_j)=\tau_j\otimes z_j,\end{eqnarray} and extended as a
$*$-algebra map, where we write $(\tau_j)=(t_1,t_1^*,t_2,t_2^*)$ for
the generators of $H=\A(\TT^2)$. This coaction makes $\A(S^7)$ into
a left $H$-comodule algebra, {\em i.e.} an algebra in the category
${}^H\mM$. It follows that $H$ also coacts on the four-sphere
algebra $\A(S^4)$ by
\begin{equation}\label{coact four sphere}
\A(S^4)\to H\otimes\A(S^4), \qquad x\mapsto 1\otimes x,\quad
\alpha \mapsto \tau_1\tau_4\otimes\alpha,\quad
\beta\mapsto\tau_2\tau_4\otimes\beta,\end{equation} making $\A(S^4)$
into an algebra in the category ${}^H\mM$ as well.
\begin{rem}
\textup{ The most general `toric' (co)action on the sphere $S^7$
would be of a four-torus. We need to restrict to $\TT^2$ in order to
have actions which are compatible with the $\SU(2)$ fibration. In
fact, we are really dealing with a double cover $\widetilde{\TT}^2
\to \TT^2$, with $\A(\widetilde{\TT}^2)$ coacting on $\A(S^7)$ and
$\A(\TT^2)$ coacting on $\A(S^4)$, as is clear from Eqs. \eqref{eqn
coact spheres} and \eqref{coact four sphere}. We shall take the
liberty of being sloppy on this point here and in the following.}
\end{rem}
The product on $\A(S^7)$ is then deformed by comodule cotwist ({\em cf.}
\S\ref{section cocycle twists}) into
$$z_j\cdot_F z_l=F(\tau_j,\tau_l)z_jz_l,\qquad
z_j\cdot_F z_l^*=F(\tau_j,\tau_l^*)z_jz_l^*.$$ Introducing the deformation
parameter
$\eta_{jl}=\mathcal{R}_F(\tau_j,\tau_l)=F^{-2}(\tau_j,\tau_l)$ given
explicitly by
\begin{equation}\label{eqn eta
matrix}(\eta_{jl})=\begin{pmatrix}1&1&\mu&\bar\mu\\1&1&\bar\mu&\mu\\\bar\mu&\mu&1&1\\\mu&\bar\mu&1&1\end{pmatrix},
\qquad \mu=e^{\ii\pi\theta},\end{equation} the deformed algebra
relations are computed to be (dropping the product symbol $\cdot_F$)
$$
z_j\, z_l=\eta_{lj}z_l\, z_j,\qquad z_j\, z_l^*=\eta_{jl}\, z_l^*\, z_j.
$$
We denote by
$\A(S^7_\theta)$ the algebra generated by
$\{z_j,z_j^*~|~j=1,\ldots,4\}$ modulo these relations. In this way,
$\A(S^7_\theta)$ is an algebra in the category ${}^{H_F}\mM$ of left
$H_F$-comodules.

Similarly, the product on $\A(S^4)$ is twisted into
$$\alpha\cdot_F\beta=F(\tau_1\tau_4,\tau_2\tau_4)\alpha\beta,
\qquad
\alpha\cdot_F\beta^*=F(\tau_1\tau_4,\tau_2^*\tau_4^*)\alpha\beta^*.$$
With deformation parameter $\lambda:=\mu^2=e^{\ii 2 \pi\theta}$, the algebra relations
become (again dropping the product symbol $\cdot_F$)
$$\alpha\beta = \lambda \beta\alpha,\quad
\alpha^*\beta^*=\lambda\beta^*\alpha^*, \quad
\beta^*\alpha=\lambda\alpha\beta^*,\quad
\beta\alpha^*=\lambda\alpha^*\beta,$$ with $x$ central. We denote by
$\A(S^4_\theta)$ the algebra generated by $\alpha$, $\beta$, $x$ and
their conjugates, subject to these relations. They make
$\A(S^4_\theta)$ into an algebra in the category ${}^{H_F}\mM$.

Since the coaction of $H$ on $\A(S^7)$ commutes with the
$\SU(2)$-action \eqref{eqn coact spheres}, the deformation of the
spheres $\A(S^7)$ and $\A(S^4)$ preserves this action and hence
there is an algebra inclusion
$\A(S^4_\theta)\hookrightarrow\A(S^7_\theta)$, giving a
noncommutative principal bundle with classical structure group
$\SU(2)$. As mentioned, the above $H$-coaction is the only one which is compatible
with the bundle structure \cite{lvs:pfns}. The elements
$\alpha$, $\beta$, $x$ and their adjoints are the entries of the
projection $\qp$  which is now given by
\begin{equation} \label{eqn basic instanton projector}
\qp:=uu^*=\tfrac{1}{2}\begin{pmatrix} 1+x & 0 & \alpha & -\bar \mu\,
\beta^* \\ 0 & 1+x & \beta & \mu\, \alpha^* \\ \alpha^* & \beta^* &
1-x & 0 \\ -\mu\, \beta & \bar\mu\, \alpha & 0 & 1-x
\end{pmatrix};
\end{equation}
note that the matrix $u$ is still of the form in Eq.~\eqref{u}.

\subsection{Noncommutative differential calculi}\label{nc diff}
There are canonical differential structures on each of the spheres
$\A(S^7_\theta)$ and $\A(S^4_\theta)$ as deformations of their
classical counterparts. They are constructed as follows.

We begin with the space $\A(\C^4)$. Let $\Omega(\C^4)$ be the usual
differential calculus on $\A(\C^4)$, generated as a commutative
differential graded algebra by the degree zero elements $z_j,z_l^*$ and
degree one elements $\D z_j, \D z_l^*$, satisfying the relations
$$
\D z_j \wedge \D z_l + \D z_l \wedge \D z_j=0,\qquad \D z_j \wedge
\D z_l^* + \D z_l^* \wedge \D z_j=0 ,
$$
for $j,l=1,\ldots,4$. The
differential $\D$ is defined by $z_j\mapsto \D z_j$, $z_l^*\mapsto
\D z_l^*$ and extended uniquely using a graded Leibniz rule. The
coaction $\Delta_L$ of Eq.~\eqref{eqn coact spheres} on $\A(\C^4)$
extends to one on the differential calculus $\Omega(\C^4)$ by defining
it to commute with the differential $\D$. We may therefore deform
the differential structure $\Omega(\C^4)$ in the same way as we did
for the algebra itself, by comodule cotwist:
\begin{multline}\label{eqn diff calc}
z_j \cdot_F \D z_l=F(\tau_j,\tau_l)z_l\D z_l, \quad
z_j\cdot_F \D z^*_l=F(\tau_j,\tau_l^*)z_j\D z^*_l,   \\
\D z_j\wedge_F \D z_l=F(\tau_j,\tau_l)\D z_j\wedge \D z_l.
\end{multline}
There is hence a
canonical differential graded algebra $\Omega(\C^4_\theta)$ for
$\A(\C^4_\theta)$, with the same generators but now subject to the
relations (again no explicit deformed product symbol)
\begin{align*}
z_j\D z_l=\eta_{lj}(\D z_l)z_j, & \qquad z_j\D z_l^*=\eta_{jl}(\D
z_l^*)z_j, \\ \D z_j \wedge \D z_l + \eta_{lj}\D z_l \wedge \D
z_j=0, & \qquad \D z_j \wedge \D z_l^*+\eta_{jl}\D z_l^* \wedge \D
z_j=0
\end{align*}
for $j,l=1,\ldots,4$. Note that the relations in the twisted
calculus $\Omega(\C^4_\theta)$ are the same as those for the
coordinate algebra $\A(\C^4_\theta)$, but now with $\D$ inserted.

Since the differential $\D$ is undeformed, the same strategy also defines a
differential calculus $\Omega(S^7_\theta)$ on $\A(S^7_\theta)$ as a
cotwist of the classical one, with the products between generators
also given by Eq.~\eqref{eqn diff calc}. Similarly, one obtains a
differential calculus $\Omega(S^4_\theta)$ on the four-sphere
$\A(S^4_\theta)$. It is generated by the degree zero elements
$\alpha$, $\alpha^*$, $\beta$, $\beta^*$, $x$ and the degree one
elements $\D\alpha$, $\D\alpha^*$, $\D\beta$, $\D\beta^*$, $\D x$, obeying relations as for the coordinate algebra $\A(S^4_\theta)$
but now with $\D$ inserted, namely
$$\alpha \, \D \beta=\lambda\D \beta \, \alpha,\qquad \beta^*
\, \D \alpha =\lambda \D \alpha \, \beta^*, \quad \D \alpha
\, \D \beta + \lambda\D\beta\, \D \alpha = 0
$$ and so on, with $x$
and $\D x$ obeying the same undeformed relations as in the classical case. The
calculus $\Omega(S^4_\theta)$ may be obtained either as the
$\SU(2)$-invariant part of $\Omega(S^7_\theta)$ or directly as a
comodule cotwist of its classical counterpart.

The torus $\TT^2$ acts on the sphere $S^4$ by isometries and hence leaves the
conformal structure invariant. As a consequence, one checks that the classical Hodge
operator (in particular on two-forms)
$*:\Omega^2(S^4)\to\Omega^2(S^4)$
is an intertwiner for the coaction $\Delta_L$ of the
torus algebra $H=\A(\TT^2)$, that is to say
$$
\Delta_L(*\omega)=(\textrm{id}\otimes *)\Delta_L(\omega), \qquad \omega \in
\Omega^2(S^4).
$$
Since the deformed differential calculus $\Omega(S^4_\theta)$ coincides as a vector space with its
undeformed counterpart $\Omega(S^4)$, we can define a Hodge operator
$*_\theta$ on $\Omega(S^4_\theta)$ by the same formula as the
classical $*$, yielding a map $*_\theta:\Omega^2(S^4_\theta) \to \Omega^2(S^4_\theta)$ which is by
construction a morphism in the category ${}^{H_F}\mM$. This is all we need when studying instantons
on $S^4_\theta$.

\section{Braided Matrix Algebras}\label{section quantum
conformal group} The previous section constructed the coordinate
algebras of the noncommutative spaces $\C^4_\theta$, $S^7_\theta$
and $S^4_\theta$ as objects in the category of left $H_F$-comodules.
In this section we observe that the various matrix algebras which
act upon these spaces may also be naturally thought of as objects in the same
category and, as a result, they are obtained using exactly the same
`quantisation' procedure.

\subsection{The classical groups $\SL(2,\HH)$ and
$\Sp(2)$} We denote by $\M(2,\HH)$ the algebra of $2\times 2$
matrices with quaternion entries; for convenience we shall write
them as $4\times 4$ matrices with complex entries. The classical
bialgebra $\A(\M(2,\HH))$ of functions on $\M(2,\HH)$ is defined to
be the commutative associative algebra generated by the coordinate
functions arranged in the following $4\times 4$ matrix
\begin{equation} \label{eqn defining
M(H)}A=\begin{pmatrix}a_{ij}&b_{ij}\\c_{ij}&d_{ij}\end{pmatrix}=\begin{pmatrix}a_1&-a_2^*&b_1&-b_2^*\\a_2&a_1^*&b_2&b_1^*\\c_1&-c_2^*&d_1&-d_2^*\\c_2&c_1^*&d_2&d_1^*\end{pmatrix}.\end{equation}
We think of this matrix as being generated by a set of
quaternion-valued functions, writing
$$a=(a_{ij})=\begin{pmatrix}a_1&-a_2^*\\a_2&a_1^*\end{pmatrix}$$ and
similarly for the other entries $b,c,d$. The $*$-structure on this
algebra is evident from the matrix (\ref{eqn defining M(H)}). We
equip $\A(\M(2,\HH))$ with the matrix bialgebra structure
$$
\Delta(A_{ij})=\sum_\alpha A_{i\alpha}\otimes A_{\alpha j}, \qquad
\ep(A_{ij})=\delta_{ij} \, \qquad \textup{for} \quad i,j=1,\ldots,4
.
$$

Of course, it is not a Hopf algebra since it does not have an
antipode (this is equivalent to saying that the matrix algebra
$\M(2,\HH)$ is not quite a group). We obtain a Hopf algebra by
passing to the quotient of $\A(\M(2,\HH))$ by the Hopf $*$-ideal
generated by the element $D-1$, where $D=\det(A)$ is the
determinant of the matrix $A$. We denote the quotient by
$\A(\SL(2,\HH))$, the algebra of functions on the group $\SL(2,\HH)$
of matrices in $\M(2,\HH)$ with determinant one. The algebra
$\A(\SL(2,\HH))$ inherits a $*$-bialgebra structure from that of
$\A(\M(2,\HH))$ and we define an antipode by
\begin{equation}\label{eqn antipode}S(A_{ij})=(-1)^{i+j}A_{ji}',$$ with $$A_{ij}':=\sum_{\sigma
\in S_3} (-1)^{|\sigma|}\ep^{\sigma_1\ldots
\sigma_{i-1}l\sigma_{i+1}\ldots \sigma_4} A_{1,\sigma_1}\ldots
A_{i-1,\sigma_{i-1}}A_{i+1,\sigma_{i+1}} \ldots
A_{4,\sigma_4},\end{equation} and $\ep^{ijkl}$ is the alternating
symbol on four elements. The notation is
$$(\sigma_1,\ldots,\sigma_{i-1},\sigma_{i+1},\ldots,\sigma_4)=\sigma(1,\ldots,l-1,l+1,\ldots,4).$$
with $\sigma$ an element of $S_3$, the permutation group on three
objects (once an index is fixed the remaining one can take only three possible values).

\begin{defn}The datum $\A(\SL(2,\HH))=(\A(\SL(2,\HH)),\Delta,\ep,S)$ constitutes a Hopf algebra.
We define the Hopf algebra
$\A(\Sp(2))$ to be the quotient of $\A(\SL(2,\HH))$ by the two-sided
$*$-Hopf ideal $\mathcal{I}$ generated by elements
\begin{equation}\label{eqn ideal gens}
\sum_\alpha(A^*)_{i\alpha}A_{\alpha j} - \delta_{ij}, \qquad
i,j=1,\ldots,4.\end{equation} In the algebra $\A(\Sp(2))$ there are
relations $A^*A=AA^*=1$, or equivalently $S(A)=A^*$.
\end{defn}

\subsection{The braided groups $\B(\SL_\theta(2,\HH))$ and $\B(\Sp_\theta(2))$}\label{section braided groups}
There is an embedding of a two-torus $\TT^2$ into $\M(2,\HH)$ as a
diagonal subgroup, given by the map
\begin{equation}\label{eqn torus embedding}\rho:\TT^2 \rightarrow \M(2,\HH), \qquad \rho(s)=\textup{diag}(e^{2\pi \ii s_1},e^{2\pi \ii s_2}),\end{equation} where $s=(e^{2\pi \ii s_1},e^{2\pi \ii s_2}) \in
\TT^2$. At the level of coordinate algebras, this inclusion becomes
a bialgebra projection
\begin{equation}\label{eqn bialg proj}\pi:\A(\M(2,\HH)) \rightarrow \A(\TT^2), \qquad
\pi(A_{ij})=\delta_{ij}\tau_j,\end{equation} where
$(\tau_j)=(t_1,t_1^*,t_2,t_2^*)$ in terms of the generators of
$H=\A(\TT^2)$. A resulting right adjoint action of $\TT^2$ on $\M(2,\HH)$ given by
$$\M(2,\HH)\times\TT^2 \to \M(2,\HH),\qquad (g,s)\mapsto
\rho(s^{-1})\cdot g \cdot\rho(s), $$
in turn dualises to the left $H$-adjoint coaction given by
\begin{equation}\label{eqn torus coaction again}\textup{Ad}_L:\A(\M(2,\HH))\rightarrow H\otimes\A(\M(2,\HH)), \qquad \textup{Ad}_L(A_{ij})=\tau_i\tau_j^*\otimes A_{ij}\end{equation} for $i,j=1,\ldots,4$ and extended as a $*$-algebra map. This coaction realises $\A(\M(2,\HH))$ as an object
in the category ${}^H\mM$ of left $H$-comodules. Since the algebra
is commutative, it follows that its product is a morphism in the
category,
$$\textup{Ad}_L(A_{ij}A_{kl})=\textup{Ad}_L(A_{ij})\textup{Ad}_L(A_{kl}).$$
The fact that the adjoint action preserves matrix multiplication in
$\M(2,\HH)$ means that the coproduct on $\A(\M(2,\HH))$ is covariant
under the coaction $\textup{Ad}_L$, {\em i.e.} $\A(\M(2,\HH))$ is an $H$-comodule
coalgebra. The same is true of the counit, whence $\A(\M(2,\HH))$ is
a bialgebra in the category ${}^H\mM$ of left $H$-comodules.
Similarly the antipodes on $\A(\SL(2,\HH))$ and $\A(\Sp(2))$
respect the $H$-coaction and so they form Hopf algebras in
the category ${}^H\mM$.

What happens to this picture under cotwisting is clear. We know from
Example~\ref{example twisted torus} that
upon twisting the torus algebra $H=\A(\TT^2)$ with a twist $F$ like the one in
\S\ref{section nc hopf fib}, as a Hopf $*$-algebra $H_F=H$, although
the coquasitriangular structure twists and the deformation takes the
form of a `quantisation functor', {\em i.e.} an isomorphism of
braided monoidal categories from ${}^H\mM$ to ${}^{H_F}\mM$. This
functor leaves objects and coactions unchanged, and hence the
adjoint coaction \eqref{eqn torus coaction again} also realises
$\A(\M(2,\HH))$ as an object in the category ${}^{H_F}\mM$. However,
we need to deform the bialgebra structures (the product and
coproduct) on $\A(\M(2,\HH))$ in order to maintain covariance.

Using Eq.~\eqref{eqn twisted prod}, the product is deformed into
\begin{equation}\label{eqn braided
prod}A_{ij}\,\underline{\cdot}\,A_{kl}=F(\tau_i\tau_j^*,\tau_k\tau_l^*)A_{ij}A_{kl},
\end{equation}
denoting the twisted product by $\underline{\cdot}$ in order to
stress its being a morphism in a braided category. We write
$\mathcal{B}(\M_\theta(2,\HH))$ for the algebra generated by the
$A_{ij}$, $i,j=1,\ldots,4$, equipped with the twisted product. Likewise, using
Eq.~\eqref{eqn twisted coprod} the coproduct is deformed on
generators into
\begin{equation}\label{eqn braided coprod}\Delta_F(A_{ij})=\sum_\alpha A_{i\alpha}\otimes A_{\alpha
j}F^{-1}(\tau_i\tau_\alpha^*,\tau_\alpha\tau_j^*)\end{equation} and
extended then as an algebra homomorphism to the braided tensor
product,
\begin{equation}\label{eqn braided
delta}\Delta_F:\B(\M_\theta(2,\HH))\to\B(\M_\theta(2,\HH))\,\underline{\otimes}\,\B(\M_\theta(2,\HH)).\end{equation}
If one defines a new set of generators of $\B(\M_\theta(2,\HH))$ by
\begin{equation}\label{redef}\widehat{A}_{ij}:=F^{-1}(\tau_i,\tau_j)A_{ij},\end{equation} then with
respect to these generators the coproduct has the standard matrix
form
$$\Delta_F(\hatA_{ij})=\sum_\alpha\hatA_{i\alpha}\otimes \hatA_{\alpha j}.$$ In
order to obtain a braided bialgebra, it is necessary to have
$\Delta_F$ respect the algebra structure of $\B(\M_\theta(2,\HH))$:
the fact that it does so is a consequence of the dual version of
\cite[Thm~2.8]{ma:qst}, although we can prove it directly as
follows. We first note that, using the Hopf bicharacter property of
$F$, one has
$$
F(\tau_i\tau_\alpha^*,\tau_\alpha\tau_j^*)=F(\tau_i,\tau_\alpha)F^{-1}(\tau_i,\tau_j)F^{-1}(\tau_\alpha,\tau_\alpha)F(\tau_\alpha,\tau_j).
$$

\begin{lem}\label{prop braided extend coproduct}The coproduct $\Delta_F$ and the product $\underline{\cdot}$ make $\B(\M_\theta(2,\HH))$ into a bialgebra in the category ${}^{H_F}\mM$ of
left $H_F$-comodules.\end{lem}

\proof By construction, the vector space $\B(\M_\theta(2,\HH))$ equipped
with the product $\underline{\cdot}$ and the coproduct $\Delta_F$ is
certainly both an algebra and a coalgebra in the category
${}^{H_F}\mM$ {\em via} the left adjoint coaction. Using the product
$\underline{\cdot}$ and the braiding in the category, we now compute
that
\begin{align*}
\Delta_F(\hatA_{ij}\,\underline{\cdot}\,\hatA_{kl})&=\Delta_F(\hatA_{ij}\hatA_{kl})F(\tau_i\tau_j^*,\tau_k\tau_l^*)\\
&=\sum_{\alpha,\beta}\hatA_{i\alpha}\hatA_{k\beta}\otimes\hatA_{\alpha j}\hatA_{\beta l}F^{-1}(\tau_i\tau_\alpha^*\tau_k\tau_\beta^*,\tau_\alpha\tau_j^*\tau_\beta\tau_l^*)F(\tau_i\tau_j^*,\tau_k\tau_l^*)\\
&=\sum_{\alpha,\beta}\hatA_{i\alpha}\, \underline{\cdot}\,\hatA_{k\beta}\otimes\hatA_{\alpha j}\,\underline{\cdot}\,\hatA_{\beta l} F^{-1}(\tau_i\tau_\alpha^*\tau_k\tau_\beta^*,\tau_\alpha\tau_j^*\tau_\beta\tau_l^*) \, \times \\
&\qquad\qquad\qquad\qquad\qquad\qquad \times \,
F(\tau_i\tau_j^*,\tau_k\tau_l^*)F^{-1}(\tau_i\tau_\alpha^*,\tau_k\tau_\beta^*)F^{-1}(\tau_\alpha\tau_j^*,\tau_\beta\tau_l^*)\\
&=\sum_{\alpha,\beta}\hatA_{i\alpha}\,\underline{\cdot}\,\hatA_{k\beta}\otimes\hatA_{\alpha
j}\,\underline{\cdot}\,\hatA_{\beta l}
F^{-2}(\tau_k\tau_\beta^*,\tau_\alpha\tau_j^*)\\&=\Delta_F(\hatA_{ij})\,\underline{\cdot}\,\Delta_F(\hatA_{kl}),
\end{align*}
where in the first equality we have used the definition of the
twisted product, in the second equality we have applied the
definition of $\Delta_F$ and in the third equality we have again
used the definition of the product. The fourth equality involves a
simplification of the terms in $F$ using its Hopf bicharacter
properties. The coproduct $\Delta_F$ is therefore an algebra map
with respect to the twisted product, whence the result.\endproof

In the same way, it follows that there are bialgebras
$\B(\SL_\theta(2,\HH))$ and $\B(\Sp_\theta(2))$ in the category
${}^{H_F}\mM$, obtained by comodule cotwist of the classical
bialgebras $\A(\SL(2,\HH))$ and $\A(\Sp(2))$. The same formula as in
Eq.~\eqref{eqn antipode}, although now using the braided product, defines an antipode for $\B(\SL_\theta(2,\HH))$ which we now
denote by $\underline{S}$. The antipode is extended as a morphism in
the category ${}^{H_F}\mM$, namely as a braided anti-algebra map
$$
\underline{S}\circ \underline{\cdot}=\underline{\cdot} \circ (\underline{S}\otimes
\underline{S})\circ \Psi,
$$
where $\Psi$ is the braiding in the category of left
$H_F$-comodules, thus making $\B(\SL_\theta(2,\HH))$ into a braided
Hopf algebra. Similarly, the formula
$$\underline{S}(A_{ij})=(A^*)_{ji},$$ extended as a braided
anti-algebra map, defines an antipode on $\B(\Sp_\theta(2))$, also
making it into a braided Hopf algebra in the category.

\section{Braided Symmetries of Noncommutative Spheres}\label{sect
braided} Now that we have constructed braided versions of the
transformation algebra $\M(2,\HH)$ and its various quotients, we are
able to show how they (co)act upon the spaces $\C^4_\theta$,
$S^7_\theta$ and $S^4_\theta$ given in \S\ref{section noncomm hopf}.
The important technical point that we illustrate is that these
coactions are themselves morphisms in the category ${}^{H_F}\mM$ and
so necessarily braided. We then construct the cobosonisations of the
transformation algebras using the procedure described
in \S\ref{section hopf algebra prelims}, the advantage
of doing so being, as mentioned, that it takes us from the realm of
braided geometry back into the realm of `ordinary' quantum groups.

\subsection{Symmetries in the braided category}
Recall that we arranged the generators of the algebra $\A(\C^4)$
into the $4\times 2$ matrix $u=(u_{ia})$ for $i=1,\dots,4$ and
$a=1,2$, given in Eq.~\eqref{u}.
In this notation the $H$-coaction can be written
$$
\A(\C^4)\to H\otimes\A(\C^4), \qquad u_{ia}\mapsto \tau_i\otimes
u_{ia}.
$$
A left coaction of the classical bialgebra $\A(\M(2,\HH))$ on
$\A(\C^4)$ is given on generators by
\begin{equation}\label{classical coaction}
\A(\C^4)\to\A(\M(2,\HH))\otimes\A(\C^4),\qquad u_{ia}\mapsto
\sum_{\beta}A_{i\beta}\otimes u_{\beta a},
\end{equation} for $i=1,\ldots,4$, $a=1,2$ and extended as a $*$-algebra map.
It is clear that this
coaction is a morphism in the category ${}^H\mM$, which becomes a
morphism in ${}^{H_F}\mM$ upon applying the quantisation functor. As
discussed in \S\ref{section cocycle twists}, the coaction itself
does not change, but we must remember that the monoidal structure is
deformed. In this way, we get a left coaction of the braided
bialgebra $\B(\M_\theta(2,\HH))$ on $\A(\C^4_\theta)$, given on
generators by
\begin{equation}\label{eqn braided coaction}
\Delta_L:\A(\C^4_\theta)\rightarrow
\B(\M_\theta(2,\HH))\,\underline{\otimes}\,\A(\C^4_\theta),\qquad
u_{ia}\mapsto \sum_{\beta}\hatA_{i \beta}\otimes u_{\beta a}.
\end{equation}
As the notation suggests, this coaction extends as an algebra
homomorphism to the braided tensor product, so that we have, for example,
\begin{align*}
u_{ia}u_{lb}&\mapsto\sum_{\beta,\gamma}(\hatA_{i\beta}\otimes
u_{\beta a})\,\underline{\cdot}\, (\hatA_{l\gamma}\otimes u_{\gamma
b})
=\sum_{\beta,\gamma}\hatA_{i\beta}\,\underline{\cdot}\,\hatA_{l\gamma}\otimes
u_{\beta a} u_{\gamma
b}\,F^{-2}(\tau_l\tau_\gamma^*,\tau_\beta),\end{align*} with
$i,l=1,\ldots,4$ and $a,b=1,2$. This argument also applies to the
braided Hopf algebra $\B(\Sp_\theta(2))$, yielding a left coaction
given by a similar expression:
\begin{equation}\label{symp coaction}
\Delta_L:\A(\C^4_\theta)\to
\B(\Sp_\theta(2))\,\underline{\otimes}\,\A(\C^4_\theta),\qquad
u_{ia}\mapsto\sum_{\beta}\hatA_{i\beta}\otimes u_{\beta a}.
\end{equation}
Now, this second coaction preserves the sphere relation of
$S^7_\theta$, since we have
\begin{align*}\sum_{\alpha}z_\alpha^*z_\alpha&\mapsto \sum_{\alpha,\beta,\gamma}\left((1\otimes
z_\beta^*)\,\underline{\cdot}\,((\hatA_{\alpha\beta})^*\otimes
1)\right)\,\underline{\cdot}\,\left(\hatA_{\alpha\gamma}\otimes
z_\gamma\right)\\&=\sum_{\alpha,\beta,\gamma}(A^*)_{\beta\alpha}\,\underline{\cdot}\,A_{\alpha\gamma}\otimes
z_{\beta}^*
z_\gamma\,F^{-2}(\tau_\alpha^*\tau_\beta,\tau_\beta^*)F^{-2}(\tau_\alpha\tau_\gamma^*,\tau_\beta^*)\\&=\sum_{\alpha,\beta,\gamma}(\hatA^*)_{\beta\alpha}\,\underline{\cdot}\,\hatA_{\alpha\gamma}\otimes
z_{\beta}^*z_\gamma\,F^{-2}(\tau_\gamma,\tau_\beta)\\&=\sum_{\beta,\gamma}\delta_{\beta\gamma}\otimes
z_\beta^*z_\gamma\,F^{-2}(\tau_\gamma,\tau_\beta)=\sum_\beta
1\otimes z_\beta^*z_\beta,
\end{align*}
and thus it descends to a coaction on the sphere
$\A(S^7_\theta)$,
\begin{equation}\label{}\Delta_L:\A(S^7_\theta)\to\B(\Sp_\theta(2))\,\underline{\otimes}\,\A(S^7_\theta),\end{equation}
defined by the same formula. Similarly, it is straightforward to
check that the $\B(\Sp_\theta(2))$-coaction restricts to the
subalgebra $\A(S^4_\theta)$ generated by the entries of the
projection $\qp=uu^*$ of Eq.~\eqref{eqn basic instanton projector}.
The entries of $\qp$ are of the form
$\qp_{kl}:=(uu^*)_{kl}=\sum_a u_{ka}(u^*)_{a l}$ for $k,l=1,\ldots,4$, so that the coaction has the
form
\begin{align}\label{braided sphere coaction}
\Delta_L&~:~\A(S^4_\theta)
\rightarrow\B(\Sp_\theta(2))\,\underline{\otimes}\,\A(S^4_\theta),
\\ \nonumber \Delta_L(\qp_{kl})&=\sum_{\beta,\gamma,a}\left(\hatA_{k\beta}\otimes
u_{\beta
a}\right)\,\underline{\cdot}\,\left((\hatA_{l\gamma})^*\otimes
(u_{\gamma a})^*\,F^{-2}(\tau_\gamma\tau_l^*,\tau_\gamma^*)\right)\\
\nonumber
&=\sum_{\beta,\gamma,a}\hatA_{k\beta}\,\underline{\cdot}\,(\hatA^*)_{\gamma
l}\otimes u_{\beta a}
(u^*)_{ a\gamma}\,F^{-2}(\tau_\gamma\tau_l^*,\tau_\gamma^*)F^{-2}(\tau_\gamma\tau_l^*,\tau_\beta)\\
\nonumber
&=\sum_{\beta,\gamma}\hatA_{k\beta}\,\underline{\cdot}\,(\hatA^*)_{\gamma
l}\otimes
\qp_{\beta\gamma}F^{-2}(\tau_\gamma\tau_l^*,\tau_\beta\tau_\gamma^*),
\end{align}
for each pair of indices $k,l=1,\ldots,4$.

\subsection{Braided conformal transformations}\label{se:bcg}
The story is similar for obtaining a coaction of
$\B(\SL_\theta(2,\HH))$ on the quantum spheres, albeit slightly more
complicated. The formula \eqref{eqn braided coaction} also defines a
left coaction $\Delta_L$ of $\B(\SL_\theta(2,\HH))$ on the algebra
$\A(\C^4_\theta)$ although, just as in the classical case, it does
not preserve the sphere relation in $\A(S^7_\theta)$. Here, the
effect of the coaction is to `inflate' the spheres, {\em i.e.} it
maps the element $r^2:=\sum_\alpha z_\alpha^*z_\alpha$ to
$$
\rho^2:=\sum_\alpha\Delta_L(z_\alpha^*z_\alpha),
$$
which is not equal to
$1\otimes \sum_\alpha z_\alpha^*z_\alpha$
in the algebra $\Delta_L(\A(\C^4_\theta))$.
Since $r^2$ is self-adjoint and central in $\A(\C^4_\theta)$, we
may evaluate it as a positive real number to give the coordinate
algebra 
of a noncommutative sphere of fixed radius $r$. As this radius
varies in $\A(\C^4_\theta)$, it sweeps out a family of
seven-spheres. Similarly, we may evaluate the central element
$\rho^2$ in $\A(\widetilde \C^4_\theta):=\Delta_L(\A(\C^4_\theta))
\subset \B(\SL_\theta(2,\HH))\,\underline{\otimes}\,\A(\C^4_\theta)$
to obtain the coordinate algebra 
of a noncommutative sphere 
of fixed radius $\rho$ : as the value of $\rho$ varies in
$\A(\widetilde \C^4_\theta)$ it sweeps out another family of
seven-spheres. The effect of the coaction $\Delta_L$ of
$\B(\SL_\theta(2,\HH))$ is to map the former family onto the latter.
We introduce the notation $\A_r(S^7_\theta):=\A(\C^4_\theta)$ and
$\A_\rho(\widetilde S^7_\theta):=\A(\widetilde \C^4_\theta)$, which
does nothing other than to stress the fact that we think of the
spaces $\C^4_\theta$ and $\widetilde \C^4_\theta$ as families of
quantum seven-spheres parameterised by the values of the functions
$r^2$ and $\rho^2$, respectively.

We find a similar picture for the four-sphere $S^4_\theta$. The
coaction of $\B(\SL_\theta(2,\HH))$ does not preserve the sphere
relation, but rather gives
$$\Delta_L(\alpha\alpha^*+\beta\beta^*+x^2)=\rho^4,$$ whence
$S^4_\theta$ is also `inflated' by the action of
$\SL_\theta(2,\HH)$.  Writing $\A_r(S^4_\theta)$ for the
$*$-subalgebra of $\A(\C^4_\theta)$ generated by $\alpha$, $\beta$,
$x$ and their conjugates, then as $r^4$ varies it sweeps out a
family of noncommutative four-spheres. Evaluating $r^4$ as a real
number yields the coordinate algebra of a noncommutative four-sphere
of radius $r^2$.

Similarly, we define $\tilde\alpha:=\Delta(\alpha)$,
$\tilde\beta:=\Delta_L(\beta)$, $\tilde x:=\Delta_L(x)$ and so
forth, writing $\A_\rho(\widetilde S^4_\theta)$ for the
$*$-subalgebra of $\A(\widetilde \C^4_\theta)$ that they generate.
As the value of $\rho^4$ varies in $\A(\widetilde\C^4_\theta)$, it
sweeps out a family of noncommutative four-spheres of radius
$\rho^2$. The effect of the coaction $\Delta_L$ is to map the former
family onto the latter.

The algebra $\A_r(S^4_\theta)$ is precisely the $\SU(2)$-invariant
subalgebra of $\A_r(S^7_\theta)$ and $\A_\rho(\widetilde
S^4_\theta)$ is the $\SU(2)$-invariant subalgebra of
$\A_\rho(\widetilde S^7_\theta)$. Consequently, there is a family of
noncommutative principal bundles parameterised by the function
$r^2$, given by the algebra inclusion
$\A_r(S^4_\theta)\hookrightarrow \A_r(S^7_\theta)$.
Similarly, there is a family of $\SU(2)$ principal bundles given by
the inclusion $\A_\rho(\widetilde S^4_\theta)\hookrightarrow
\A_\rho(\widetilde S^7_\theta)$. The discussion above shows that the coaction of
$\B(\SL_\theta(2,\HH))$ serves to map the former family onto the
latter. For further details, we refer to \cite{lprs:ncfi,bl:adhm}.

\begin{prop}\label{conf trans}
With $*_\theta:\Omega^2(S^4_{\theta;r^2})\to
\Omega^2(S^4_{\theta;r^2})$ the Hodge operator on the sphere
$S^4_{\theta;r^2}$ of fixed radius $r^2$, the braided Hopf algebra
$\B(\SL_\theta(2,\HH))$ coacts on $\Omega^2(S^4_{\theta;r^2})$ by
conformal transformations, that is
$$
\Delta_L(*_\theta\,\omega)=(\id\otimes
*_\theta)\Delta_L(\omega) \qquad \text{for all} \quad \omega\in\Omega^2(S^4_{\theta;r^2}).
$$
\end{prop}

\proof The coaction $\Delta_L$ of $\B(\SL_\theta(2,\HH))$ is
extended to forms $\Omega(S^4_{\theta;r^2})$ by requiring it to
commute with $\D$, namely
$\Delta_L(\D\omega)=(\id\otimes\D)\Delta_L(\omega)$ for all
$\omega\in\Omega(S^4_{\theta;r^2})$. Now the coaction $\Delta_L$ is
given by the classical action of $\SL(2,\HH)$ on
$\Omega(S^4_{\theta=0;r^2})$ as vector spaces and only the products
on $\A(\SL(2,\HH))$ and $\A(S^4)$ are deformed. Since $*_\theta$
coincides with the classical Hodge operator $*$ on
$\Omega(S^4_{\theta;r^2})\simeq \Omega(S^4_{\theta=0;r^2})$ as
vector spaces, the result follows from the classical fact that
$\SL(2,\HH)$ acts on $S^4$ by conformal transformations.\endproof

In order to proceed we need to slightly enlarge all of our algebras
and assume that the quantity $r^2$ is invertible with inverse
element $r^{-2}$, and shall henceforth assume that this is done without change of notation.
In terms of our families of seven-spheres, this means that we now
think of $\A_r(S^7_\theta)$ in the same way as before but without
the `origin' in $\C^4_\theta$, which corresponds to the `sphere of
radius zero'.

We also define an inverse $\rho^{-2}$ for the quantity $\rho^2$,
extending the coaction to the new elements by
$\rho^{-2}:=\Delta_L(r^{-2})$. This gives a well-defined coaction,
$$
\Delta_L:\A_r(S^7_\theta)\to\B(\SL_\theta(2,\HH))\,\underline{\otimes}\,\A_r(S^7_\theta),
$$
for which $\A_r(S^7_\theta)$ is a braided
$\B(\SL_\theta(2,\HH))$-comodule algebra. Following the above, for
the image under $\Delta_L$ we write $\A_\rho(\widetilde
S^7_\theta):=\Delta_L(\A_r(S^7_\theta))$, noting that $\rho^2$ and
$\rho^{-2}$ are central in $\A_\rho(\widetilde S^7_\theta)$ since
$r^2$ and $r^{-2}$ are central in $\A_r(S^7_\theta)$. In these new
terms, the construction of the defining projector of
$\A_r(S^4_\theta)$ needs only a minor modification. We now use
\begin{equation}\label{rescaled u}\sfu:=\begin{pmatrix}z_1&z_2&z_3&z_4\\-z_2^*&z_1^*&-z^*_4&z_3^*\end{pmatrix}^{\textrm{t}},
\end{equation}
to denote the same matrix as in Eq.~\eqref{u} but now without
imposing the sphere relations, so that now we have $\sfu^*\sfu=r^2$.
It follows that the matrix
\begin{equation}\label{rescaled p}\qp:=\sfu \, r^{-2}\sfu^*=\tfrac{1}{2}r^{-2}\begin{pmatrix} r^2+x & 0 &
\alpha
& -\bar \mu \,\beta^* \\ 0 & r^2+x & \beta & \mu \, \alpha^* \\
\alpha^* & \beta^* & r^2-x & 0 \\ -\mu \,\beta & \bar\mu \, \alpha &
0 & r^2-x\end{pmatrix}\end{equation} is a projection whose entries
generate the $\SU(2)$-invariant subalgebra of $\A_r(S^7_\theta)$.
Moreover, there is a well-defined left coaction of
$\B(\SL_\theta(2,\HH))$ on the algebra $\A_r(S^4_\theta)$ generated
by the entries of this matrix, given by the same formula as in
Eq.~\eqref{braided sphere coaction}. Writing $\widetilde
\sfu:=\Delta_L(\sfu)$, the image of the projector $\qp$ under the
braided coaction $\Delta_L$ is computed to be
\begin{equation} \label{eqn coacted instanton projector}
\widetilde{\qp}:=\widetilde\sfu
\rho^{-2}\widetilde\sfu^*=\tfrac{1}{2}\rho^{-2}\begin{pmatrix}
\rho^2+\tilde x & 0 &
\tilde \alpha & -\bar \mu \,\tilde \beta^* \\ 0 & \rho^2+x & \tilde \beta & \mu \, \tilde \alpha^* \\
\tilde \alpha^* & \tilde \beta^* & \rho^2-x & 0 \\ -\mu \,\tilde
\beta & \bar\mu \, \tilde \alpha & 0 & \rho^2-\tilde x
\end{pmatrix}\end{equation} and defines an element in the K-theory
of the algebra
$\B(\SL_\theta(2,\HH))\,\underline{\otimes}\,\A_r(S^4_\theta)$. By
construction we have that $\widetilde \qp=\Delta_L(\qp)$, with the
entries of $\widetilde \qp$ explicitly given by
$$
\widetilde\qp_{kl}=\sum_{\alpha,\beta}\rho^{-2}\hatA_{k\alpha}\,\underline{\cdot}\,(\hatA^*)_{\beta
l}\otimes(\sfu\sfu^*)_{\alpha\beta}\,F^{-2}(\tau_\beta\tau_l^*,\tau_\alpha\tau_\beta^*),
$$ whence we shall think of $\widetilde \qp$ as a `braided
family' of projections parameterised by the algebra
$\B(\SL_\theta(2,\HH))$.

\subsection{The cobosonised transformation algebra}
As mentioned, when working in the braided category ${}^{H_F}\mM$, we have to remember not only
that the algebras are twisted, but that the coactions (as braided
morphisms in the category) are also twisted since they involve the
tensor product structure of the category. This can be
computationally rather awkward, so it is useful to remember that the
braided left comodules for $\B(\SL_\theta(2,\HH))$ (similarly for
$\B(\Sp_\theta(2))$) are in one-to-one correspondence with left
comodules for its cobosonisation, which is a Hopf algebra in the
`ordinary' sense that we now compute.

From \S\ref{section braided groups} the left adjoint
coaction of $H_F$ on $\B(\M_\theta(2,\HH))$ is given by
\begin{equation}\label{eqn recall torus coaction}\B(\M_\theta(2,\HH))\rightarrow H_F \otimes\B(\M_\theta(2,\HH)), \qquad \widehat A_{ij} \mapsto \tau_i\tau_j^*\otimes \hatA_{ij},
\end{equation} for $i,j =1,\ldots,4$. This coaction
makes $\B(\M_\theta(2,\HH))$ into a coalgebra in the category
${}^{H_F}\mM$ of left $H_F$-comodules, which means that we may
construct the associated crossed coproduct coalgebra
$\B(\M_\theta(2,\HH))\lcocross H_F$ defined in \S\ref{section
hopf algebra prelims}. As a vector space it is just
$\B(\M_\theta(2,\HH))\otimes H_F$, with the cross coproduct defined
in Eq.~\eqref{eqn cross coproduct} working out to be
\begin{equation}\label{eqn explicit cross coproduct formula}\Delta(\widehat A_{ij}\otimes h)=\sum_\alpha\widehat A_{i\alpha} \otimes
\tau_\alpha\tau_j^*h \otimes\widehat A_{\alpha j}\otimes
h\end{equation} on group-like elements $h\in H_F$
and extended by linearity.

With the coquasitriangular structure $\mathcal{R}_F=F^{-2}$ on
$H_F$, we have also an $H_F$-action,
\begin{equation}\label{eqn-torus-action}
H_F \times \B(\M_\theta(2,\HH)) \rightarrow \B(\M_\theta(2,\HH)) ,
\qquad h\tr\widehat A_{ij}=\widehat A_{ij}
F^{-2}(\tau_i\tau_j^*,h),
\end{equation}
with $i,j =1,\ldots,4$. Thus, we may canonically view $\B(\M_\theta(2,\HH))$ as an object
in the category ${}^{H_F}_{H_F}\mC$ of crossed $H_F$-modules. It follows that we may
construct the associated cross product algebra
$\B(\M_\theta(2,\HH))\lcross H_F$, which also has
$\B(\M_\theta(2,\HH))\otimes H_F$ as an underlying vector space. The
cross product defined in Eq.~\eqref{eqn cross product} works out as
$$(\widehat A_{ij}\otimes g)(\widehat A_{kl}\otimes h)=\widehat A_{ij}\,\underline{\cdot}\,(g\tr\widehat A_{kl})\otimes gh=\widehat
A_{ij}~\underline{\cdot}~\widehat A_{kl}\otimes
gh\,F^{-2}(\tau_k\tau_l^*,g)$$ for group-like elements $g,h,\in
H_F$. From \S\ref{section hopf algebra prelims}, we know that we
have constructed the cobosonised bialgebra
$\B(\M_\theta(2,\HH))\lbiprod H_F$.

\begin{rem}\label{decompose}
\textup{From the general theory of biproducts,
$\B(\M_\theta(2,\HH))\lbiprod H_F$ contains $H_F$ as a sub-Hopf
algebra {\em via} the projection $\pi_H:=\underline{\ep}\otimes\id$.
The subalgebra $\B(\M_\theta(2,\HH))$ is recovered as the algebra of
coinvariants under the right coaction $(\id\otimes\pi_H)\circ
\Delta$. Moreover, it is interesting to note that the cobosonisation
$\B(\M_\theta(2,\HH))\lbiprod H_F$ contains the transformation
bialgebra $\A(\M_\theta(2,\HH))$ constructed in \cite{lprs:ncfi} (in
fact it is isomorphic to the double cross product
$\A(\M_\theta(2,\HH))\dcross H_F$, {\em cf}. \cite{ma:qd,ma:book}).
This is to be expected, since $\A(\M_\theta(2,\HH))$ was constructed
as the {\em universal} transformation bialgebra of
$\A(\C^4_\theta)$.}
\end{rem}

Similarly, we may construct the cobosonisations
$\B(\SL_\theta(2,\HH))\lbiprod H_F$ and $\B(\Sp_\theta(2))\lbiprod
H_F$ of the braided Hopf algebras $\B(\SL_\theta(2,\HH))$ and
$\B(\Sp_\theta(2))$. The antipodes on these Hopf algebras are given
by Eq.~\eqref{eqn biproduct antipode} and come out on generators to
be
$$S(\widehat A_{ij}\otimes h)=(1\otimes \tau_j\tau_i^*h^*))(\underline{S}(\hatA_{ij})\otimes 1),$$ with $\underline{S}$
the braided antipode of $\B(\SL_\theta(2,\HH))$ and of
$\B(\Sp_\theta(2))$.

As a result of this construction, there is a coaction of
$\B(\M_\theta(2,\HH))\lbiprod H_F$ on $\A_r(S^7_\theta)$ and its
various subalgebras. Explicitly, we have a coaction
\begin{equation}\label{cob coact}\A_r(S^7_\theta)\to
(\B(\SL_\theta(2,\HH))\lbiprod H_F)\otimes\A_r(S^7_\theta),\qquad
\sfu_{ia}\mapsto \sum_\beta\hatA_{i\beta}\otimes \tau_\beta\otimes
\sfu_{\beta a},\end{equation} on $\A_r(S^7_\theta)$, obtained as the
coaction of $\B(\M_\theta(2,\HH))$ followed by the coaction of
$H_F$.
Just as we did for $\B(\SL_\theta(2,\HH))$ in Eq.~\eqref{eqn coacted
instanton projector}, we check that this descends to a well-defined
coaction on the family of four-spheres $\A_r(S^4_\theta)$. Indeed,
coacting upon the projection $\qp$ in this way yields a projection,
denoted $\widetilde \Qp$, with entries in the algebra
$(\B(\SL_\theta(2,\HH))\lbiprod H_F)\otimes \A_r(S^4_\theta)$. These
entries are given explicitly by
\begin{equation}\label{cobos proj}\widetilde \Qp_{kl}=\sum_{\alpha,\beta}\rho^{-2}\hatA_{k\alpha}\,\underline{\cdot}\,(\hatA^*)_{\beta
l}\otimes\tau_\alpha\tau_\beta^*\otimes(\sfu\sfu^*)_{\alpha\beta}\,F^{-2}(\tau_\beta\tau_l^*,\tau_\alpha\tau_\beta^*).\end{equation}
We think of $\widetilde \Qp$ as a noncommutative family of
projections parameterised by the algebra
$\B(\SL_\theta(2,\HH))\lbiprod H_F$.

\section{Noncommutative Families of Instantons}
\label{section gauge theory} We are now in a position to apply the
abstract theory described in previous sections to the subject of
interest: the construction of instanton connections on the sphere
$S^4_\theta$. We begin by recalling the theory of anti-self-dual
connections on $S^4_\theta$ and what it means for two such
connections to be gauge equivalent. We then extend this by
discussing the notion of equivalent families of connections over
$S^4_\theta$ and, in particular, of families of instantons.

We then provide some examples of families of instanton connections.
The first example comes from a basic instanton: by acting upon this
with various symmetry groups (namely the torus $H_F=\A(\TT^2)$, the
braided groups $\B(\SL_\theta(2,\HH))$ and $\B(\Sp_\theta(2))$, as
well as their cobosonisations $\B(\SL_\theta(2,\HH))\lbiprod H_F$
and $\B(\Sp_\theta(2))\lbiprod H_F$) we generate a variety of
different families and discuss which of them are equivalent.

\subsection{Connections and gauge equivalence}\label{gauge freedom}
Here we briefly recall the notion of gauge equivalence for
connections on vector bundles over the four-sphere $S^4_\theta$, the latter
equipped with the differential calculus
$(\Omega(S^4_\theta), \D)$ defined in \S\ref{nc diff}.

Let $\E$ be a finitely-generated projective right
$\A(S^4_\theta)$-module endowed with an $\A(S^4_\theta)$-valued
Hermitian structure $\la\,\cdot\,|\,\cdot\,\ra$; this is assumed to
be self-dual, meaning that every right $\A(S^4_\theta)$-module homomorphism
$\phi:\E\to \A(S^4_\theta)$ is represented by some element $\eta\in\E$ under the
assignment $\phi(\,\cdot\,)=\la \eta|\,\cdot\,\ra$. A connection on
$\E$ is a linear map $\n:\E\to \E\otimes_{\A(S^4_\theta)}
\Omega^1(S^4_\theta)$ satisfying the Leibniz rule
$$\n(\xi x)=(\n \xi)x + \xi \otimes \D x \qquad \text{for all} ~~\xi
\in \E, ~x \in \A(S^4_\theta).$$ The connection $\n$ is said to be
compatible with the Hermitian structure on $\E$ if it obeys
$$
\la \n \xi|\eta\ra + \la\xi|\n \eta \ra= \D \la \xi|\eta \ra \qquad
\text{for all} ~~\xi, \eta \in \E, ~x \in \A(S^4_\theta).
$$ By assumption,
$\E$ is a direct summand of a free module, that is
$\E=P(\C^N\otimes \A(S^4_\theta))$ for some $P\in
\End_{\A(S^4_\theta)}(\E)$, $P^2=P=P^*$, which we use to define the
so-called Grassmann connection $\n_0:=P\circ \D$ on $\E$. It is
straightforward to check that $\n_0$ is a compatible connection. Any
other compatible connection on $\E$ is of the form $\n=\n_0+\alpha$,
where $\alpha$ is a skew-adjoint element of
$\textup{Hom}_{\A(S^4_\theta)}(\E,\E\otimes_{\A(S^4_\theta)}
\Omega^1(S^4_\theta))$.

The curvature of $\n$ is the $\End_{\A(S^4_\theta)}(\E)$-valued
two-form $\n^2$, which in the case of the Grassmann connection
$\n_0$ is easily computed to be $\n_0^2=P(\D P)^2$. More generally,
the curvature of $\n=\n_0+\alpha$ comes out to be
$$\n^2=P(\D P)^2+P(\D \alpha)P + \alpha^2.$$ The curvature $\n^2$ is
said to be {\em anti-self-dual} if it satisfies the equation
$$*_\theta\n^2=-\n^2,$$ where
$*_\theta:\Omega^2(S^4_\theta)\to\Omega^2(S^4_\theta)$ is the Hodge
operator on two-forms; if this is the case we say that the connection $\n$ is an {\em instanton}.

The gauge group of $\E$ is defined to be
$$
\mathcal{U}(\E):=\left\{ U \in \textup{End}_{\A(S^4_\theta)}(\E)~|~ \la
U\xi|U\eta\ra=\la\xi|\eta\ra ,~\text{for all}~ \xi,\eta \in \E
\right\}.
$$
The gauge group $\mathcal{U}(\E)$ acts upon the space of compatible
connections by $$\n \mapsto \n^{U}:=U\n U^*$$ for each compatible
connection $\n$ and each element $U$ of $\mathcal{U}(\E)$. Of
course, $\n^U$ is not a `new' connection, rather it expresses $\n$
on  the `transformed bundle' $U\E$. Thus we say that a pair of
connections $\n_1$, $\n_2$ on $\E$ are {\em gauge equivalent} if
they are related by such a gauge transformation $U$. In terms of the
decomposition $\n=\n_0+\alpha$, one finds $\n^U= \n_0+\alpha^U$,
where $\alpha^U:=U (\n_0 U^*) + U\alpha U^*$. The curvature
transforms to $(\n^U)^2= U\n^2 U^*$, so that if $\n$ is an
instanton connection then so is $\n^U$.

\subsection{Families of instantons} Let $A$ be a unital $*$-algebra over $\C$.
In this section we shall investigate what it means to have a family
of connections parameterised by the algebra $A$ and define when two
such families are equivalent. We begin with the notion of a family
of vector bundles parameterised by an algebra.

\begin{defn}\label{fam of bundles}A {\em family of Hermitian vector bundles} over $S^4_\theta$ parameterised by the algebra $A$ is
a finitely-generated projective right module $\E$ over the algebra
$A\otimes\A(S^4_\theta)$ equipped with an
$A\otimes\A(S^4_\theta)$-valued Hermitian structure
$\la\,\cdot\,|\,\cdot\,\ra$. We shall give any such a module {\em
via} a self-adjoint idempotent $\Pp\in \textup{M}_N(A\otimes
\A(S^4_\theta))$, $\Pp^2=\Pp=\Pp^*$, for which $\E:=\Pp(A\otimes
\A(S^4_\theta))^N$.\end{defn}

We write $A\otimes\Omega^1(S^4_\theta)$ for the tensor product
bimodule over the algebra $A\otimes \A(S^4_\theta)$ and extend the
differential $\D$ on $\A(S^4_\theta)$ to $A\otimes \A(S^4_\theta)$
as $\id\otimes \D$.

\begin{defn}\label{fam of cons}A {\em family of connections} over $S^4_\theta$
parameterised by the algebra $A$ consists of a family of Hermitian
vector bundles $\E:=\Pp(A\otimes \A(S^4_\theta))^N$ over
$S^4_\theta$, together with a linear map
$$\n:\E\to \E\otimes_{A\otimes \A(S^4_\theta)} (A\otimes\Omega^1(S^4_\theta))
\simeq
 \E \otimes_{\A(S^4_\theta)} \Omega^1(S^4_\theta) $$
obeying the Leibniz rule $$\n(\xi x)=(\n \xi)x + \xi \otimes (\id\otimes\D) x$$ for all
$\xi \in \E$, $x \in A\otimes \A(S^4_\theta)$. The family is said to
be {\em compatible} with the Hermitian structure if it obeys $\la \n
\xi|\eta\ra + \la\xi|\n \eta \ra= (\id\otimes\D) \la \xi|\eta \ra$
for all $\xi \in \E$, $x \in A\otimes \A(S^4_\theta)$.\end{defn}

On the family of Hermitian vector bundles
$\E=\Pp(A\otimes \A(S^4_\theta))^N$ over $S^4_\theta$, there is the associated {\em
family of Grassmann connections} $\n_0=\Pp\circ (\textup{id}\otimes
\D)$. More generally, we can always express a family of connections
in the form $\n=\n_0+\omega$, where $\omega$ is an element of
$\End_{A\otimes\A(S^4_\theta)} (\E, \E\otimes_{A\otimes
\A(S^4_\theta)} (A\otimes\Omega^1(S^4_\theta))) \simeq  \End_{A\otimes\A(S^4_\theta)} (\E, \E \otimes_{\A(S^4_\theta)} \Omega^1(S^4_\theta))$.

\begin{defn}\label{fam inst}
A {\em family of instantons} over $S^4_\theta$ is a family of compatible
connections $\n$ over $S^4_\theta$ whose curvature $\n^2$ obeys the
anti-self-duality equation
$$(\textup{id}\otimes *_\theta)\n^2 = -\n^2,$$ where $*_\theta$ is the
Hodge operator on $\Omega^2(S^4_\theta)$.
\end{defn}

We also need to generalise gauge equivalence to incorporate families
of connections.

\begin{defn}\label{equiv fam}
Let $\E:=\Pp(A\otimes\A(S^4_\theta))^N$ be a family of Hermitian vector bundles parameterised
by the algebra $A$. The {\em gauge group} of $\E$ is
$$\mathcal{U}(\E):=\{U\in \End_{A\otimes\A(S^4_\theta)}(\E)~|~\la
U\xi|U\eta\ra = \la\xi|\eta\ra\, ~\text{for all}~\xi,\eta\in\E\}.$$ We
say that two families of compatible connections $\n_1$, $\n_2$ on $\E$
are {\em equivalent families} and write $\n_1\sim \n_2$ if they are
related by the action of the gauge group, {\em i.e.} there exists
$U\in \mathcal{U}(\E)$ such that $\n_2=U\n_1U^*$.\end{defn}

\begin{rem}
\textup{Where $A=\C$ ({\em i.e.} for a family parameterised by a one-point space), the above relation reduces to the usual definition of
gauge equivalence of connections. In the case where the families
$\n_1$, $\n_2$ are Grassmann families associated to projections
$P_1,P_2\in\M_N(A\otimes\A(S^4_\theta))$, equivalence means that
$P_2=UP_1U^*$ for some unitary $U$.
}
\end{rem}

\begin{prop}\label{mod func}Let $\n$ be a family of connections over $S^4_\theta$ parameterised by the algebra $A$. Then for each unital $*$-algebra morphism
$\phi:A\to B$ there is an induced family $\phi_*\n$ of connections
parameterised by the algebra $B$. This operation obeys the
functorial properties
$$(\phi_1 \circ \phi_2)_*=\phi_1{}_*\circ\phi_2{}_*, \qquad
(\id_A)_*=\textup{id},$$ and is compatible with the gauge equivalence $\sim$ in
that $\n_1\sim\n_2$ implies $\phi_*\n_1\sim \phi_*\n_2$.\end{prop}

\proof Let $\E_A$ be a finitely-generated projective
$A\otimes\A(S^4_\theta)$-module and let $\n$ be a connection on
$\E_A$ as in Definition~\ref{fam of cons}. If $\phi:A\to B$ is a
unital $*$-algebra map then we define a finitely-generated
projective right $B\otimes\A(S^4_\theta)$-module $\E_B$ by
$$\E_B:=\E_A\otimes_{A\otimes\A(S^4_\theta)}(B\otimes
\A(S^4_\theta)),$$ where $B\otimes \A(S^4_\theta)$ is thought of as
a left $A\otimes \A(S^4_\theta)$-module {\em via} the map
$\phi\otimes\id_{\A(S^4_\theta)}$. The induced connection $\phi_*\n$
is defined by
$$\phi_*\n:=\n\otimes \id_B\otimes\id_{\A(S^4_\theta)}+
\id_{\E_A}\otimes\id_B\otimes \D$$ with respect to the above
decomposition of $\E_B$. The functorial properties of $\phi_*$ are
obvious; the fact that $\phi_*$ respects unitary equivalence is also
clear.\endproof

\begin{prop}\label{inj curv}Let $\n$ be a family of connections over $S^4_\theta$
parameterised by the algebra $A$ and let $\phi:A\to B$ be a morphism of
$*$-algebras. Then the curvature of the induced family $\phi_*\n$ is
equal to the curvature of $\n$. In particular, if $\n$ is a family
of instantons then so is $\phi_*\n$.\end{prop}

\proof The curvature of $\phi_*\n$ is computed as follows. For each
$\xi\in \E_A$, $b\in B$, $y\in\A(S^4_\theta)$ we have
\begin{align*}
(\phi_*\n)^2(\xi\otimes b\otimes y)&=(\phi_*\n)\left((\n \xi)\otimes b\otimes y
+ \xi\otimes b\otimes \D y \right)\\
&=(\n^2 \xi)\otimes b\otimes y -(\n \xi)\otimes b\otimes \D y +
 (\n \xi)\otimes b\otimes \D y + \xi\otimes b\otimes \D^2 y\\
 &=(\n^2 \xi)\otimes b\otimes y,\end{align*}
 where in
the second line we have extended $\phi_*\n$ using a graded Leibniz
rule, as required for it to be well-defined on one-forms. This
simply says that, since the curvature $\n^2$ is
$A\otimes\A(S^4_\theta)$-linear, it is unaffected when we extend the
scalars to $B\otimes\A(S^4_\theta)$. It follows that if $\n^2$ is
anti-self-dual, then so is $(\phi_*\n)^2$.\endproof

\begin{rem}\label{def mod}
\textup{Writing $\mathsf{Alg}$ for the category of unital
$*$-algebras over $\C$ and $\mathsf{Set}$ for the category of sets,
these two propositions say that we have a covariant functor
$\mathcal{F}: \mathsf{Alg}\to \mathsf{Set}$. The functor maps each
algebra $A$ to the set $\mathcal{F}(A)$ of equivalence classes of
families of instantons parameterised by $A$. The functor
$\mathcal{F}$ is called the {\em moduli functor}. If we so wish, we
may restrict this functor to the sub-category ${}^{H_F}\mathsf{Alg}$
consisting of unital left $H_F$-comodule $*$-algebras ({\em cf}.
Appendix~\ref{families}).}
\end{rem}

Observe that Definition~\ref{equiv fam} only defines an equivalence
relation on families parameterised by the same algebra $A$, whereas
Proposition~\ref{mod func} provides us with a notion of equivalence for
families of instantons which are parameterised by {\em different}
algebras. Indeed, if $\n_1$ and $\n_2$ are families parameterised by
algebras $A_1$, $A_2$, we can think of them as being equivalent
if there exists an algebra $B$ and a pair of morphisms
$\phi_1:A_1\to B$, $\phi_2:A_2\to B$ such that $\phi_1{}_*\n_1\sim
\phi_2{}_*\n_2$.

\subsection{The basic instanton}\label{basic eq}
We now turn to the explicit construction of families of
connections on $S^4_\theta$, beginning with a review of the basic
instanton constructed in \cite{lvs:pfns}. With $\qp$ the projection
in Eq.~\eqref{rescaled p}, we consider the vector bundle and
Grassmann connection associated to the complementary projection
$\pp:=1-\qp$.
The Grassmann connection
\begin{equation}\label{eqn basic}
\n=\pp\circ \D:\E \rightarrow
\E\otimes_{\A(S^4_\theta)}\Omega^1(S^4_\theta)
\end{equation} has
curvature $\n^2=\pp(\D \pp)^2$ which is known to be anti-self-dual,
$$*_\theta (\pp(\D \pp)^2)=-\pp(\D \pp)^2,$$ and is hence an
instanton which we call the {\em basic instanton}. Using
noncommutative index theory, the topological charge of the
projection $\pp$ is computed to be equal to $-1$. The reason for
going from the projection $\qp$ to $\pp$ is to obtain a connection
with anti-self-dual curvature on a bundle with a fixed rank, just 2
for the case studied. In fact the Grassmann connection defined by
$\qp$ is known to have self-dual curvature, and would then qualify
to be an anti-instanton, given a happy coincidence that is unique to
the case of the lowest value of the instanton charge, coming from
the fact that the bundles described by the complementary projections
$\qp$ and $\pp$ happen to have the same rank. As we shall see when
we come to consider instantons of higher topological charge, the
crucial component in this construction is indeed the use of
complementary projections to obtain instantons, but the
corresponding bundles do not have equal rank.

The simplest way to generate new connections is to act upon the
sphere $S^4_\theta$ by a group of symmetries and look at what
happens to the basic instanton as a result. The first example of
such a symmetry group that we encountered was the two-torus $\TT^2$,
whose action was encoded in \S\ref{section nc hopf fib} as a
coaction
\begin{equation}\label{hcoas4}
\Delta_L:\A(S^4_\theta)\to H_F\otimes \A(S^4_\theta),\qquad
\qp_{kl}\mapsto \tau_k\tau_l^*\otimes \qp_{kl},
\end{equation}
for
$k,l=1,\ldots,4$. As mentioned above, we are more interested in the
complementary projection $\pp$, which transforms in the same way
under the coaction of $H_F$:
$$\Delta_L(\pp_{kl})=\Delta_L(\delta_{kl}-\qp_{kl})=\tau_k\tau_l^*\otimes(\delta_{kl}-\qp_{kl})=\tau_k\tau_l^*\otimes \pp_{kl}.$$ This leads immediately to the
following fact.

\begin{prop}\label{torus gauge}With $\Delta_L$ the coaction of $H_F$ on $\A(S^4_\theta)$ given in \eqref{hcoas4}, the projection
$\pp':=\Delta_L(\pp)$ defines a family $\n':=\pp'\circ
(\id\otimes\D)$ of instantons parameterised by the algebra $H_F$.
The family $\n'$ is equivalent to the basic instanton $\n:=\pp\circ
\D$.\end{prop}

\proof It is not difficult to check that $\pp':=\Delta_L(\pp)$ is a
projection with entries in the algebra $H_F\otimes\A(S^4_\theta)$.
We may view $\pp$ as a projection in
$\M_4(H_F\otimes\A(S^4_\theta))$ using the algebra map
$\A(S^4_\theta)\hookrightarrow H_F\otimes \A(S^4_\theta)$ defined by
$a\mapsto 1\otimes a$ for each $a \in \A(S^4_\theta)$. Taking $U$ to
be the unitary matrix
\begin{equation}\label{unitary}U:=\textup{diag}~(\tau_1\otimes 1,\tau_2\otimes 1,\tau_3\otimes
1,\tau_4\otimes 1)\in\M_4(H_F\otimes\A(S^4_\theta)) ,
\end{equation}
with $\tau_j$ the generators of $H_F$, it is straightforward to
verify that $\pp'=U(1\otimes\pp)U^*$ and so the two families are
equivalent. It follows immediately that the family $\n'$ also has
anti-self-dual curvature.\endproof

We immediately see how to generate other families of instantons
which are equivalent to the basic one. We are not limited to
conjugating $\pp$ simply by generators of $H_F$ as in
Eq.~\eqref{unitary}: more generally we can take any quadruple
$(u_1,u_2,u_3,u_4)$ of unitary elements in $H_F$ with $u_1=u_2^*$,
$u_3=u_4^*$ and set
$$U=\textup{diag}~\left(u_1\otimes 1,u_2\otimes 1,u_3\otimes
1,u_4\otimes 1\right)\in\M_4(H_F\otimes\A(S^4_\theta)).$$ The
resulting projection $U(1\otimes\pp)U^*$ is by definition equivalent
to $\pp$. We still get a family of instantons parameterised by
$H_F$, although it is a different parameterisation from the one in
Proposition~\ref{torus gauge}. Since the topological charge of the
projection $\pp$ is invariant under unitary equivalence, these
families also have charge equal to $-1$.

\subsection{A noncommutative family of instantons}\label{nc family}Next we examine the effect of the coaction
$\Delta_L:\A_r(S^4_\theta)\to \left(\B(\SL_\theta(2,\HH))\lbiprod
H_F\right)\otimes\A_r(S^4_\theta)$ of the cobosonised transformation
algebra on the basic instanton, once again stressing that for a
well-defined coaction we have to work not with $\A(S^4_\theta)$ but
with the entire family of spheres $\A_r(S^4_\theta)$. Recall from
Eq.~\eqref{cob coact} that the coaction is given by
$$\Delta_L:\A_r(S^4_\theta)\to\left(\B(\SL_\theta(2,\HH))\lbiprod
H_F\right)\otimes\A_r(S^4_\theta)$$ on $\A_r(S^4_\theta)$ and hence
on the projection $\qp$, yielding the family of projections
$\widetilde \Qp$ described in Eq.~\eqref{cobos proj}:
\begin{equation}\label{qtproj}
\widetilde \Qp_{kl}=\sum_{\alpha,\beta}\rho^{-2}\hatA_{k\alpha}\,\underline{\cdot}\,(\hatA^*)_{\beta
l}\otimes\tau_\alpha\tau_\beta^*\otimes
(\sfu\sfu^*)_{\alpha\beta}\,F^{-2}(\tau_\beta\tau_l^*,\tau_\alpha\tau_\beta^*).\end{equation}
There is a similar coaction of the same Hopf algebra on the
projection $\pp$, which may be expressed by writing $\pp=1-\qp$ and
computing
$$
\Delta_L(\pp_{kl})=\Delta_L(\delta_{kl}-\qp_{kl})=1\otimes \delta_{kl}-\widetilde \Qp_{kl}.
$$
We denote the resulting projection by $\widetilde
\Pp:=\Delta_L(\pp)$. Extending the exterior derivative to
$(\B(\SL_\theta(2,\HH))\lbiprod H_F)\otimes\A_r(S^4_\theta)$ as
$\id\otimes \D$, we get the following result.

\begin{prop}\label{big fam}
The family $\widetilde \n$ of Grassmann connections
defined by $\widetilde \n:=\widetilde \Pp \circ (\id\otimes\D)$ has
anti-self-dual curvature, that is
$$(\textup{id}\otimes
*_\theta)\widetilde \Pp((\textup{id}\otimes \D)\widetilde \Pp)^2=-\widetilde \Pp((\textup{id}\otimes \D)\widetilde
\Pp)^2.
$$
\end{prop}

\proof By definition, the coaction of $\B(\SL_\theta(2,\HH))\lbiprod
H_F$ on $\A_r(S^4_\theta)$ is given by first coacting with
$\B(\SL_\theta(2,\HH))$ followed by coacting with $H_F$. From
Proposition~\ref{conf trans}, the braided group
$\B(\SL_\theta(2,\HH))$ coacts by conformal transformations and so
commutes with the Hodge structure $*_\theta$, hence preserving the
anti-self-duality. By Proposition~\ref{torus gauge}, the Hopf
algebra $H_F$ coacts unitarily and hence preserves the curvature,
and we know that the $\B(\SL_\theta(2,\HH))$-coaction intertwines
the $H_F$-coaction. Putting these coactions together gives a family
of instantons parameterised by the algebra
$\B(\SL_\theta(2,\HH))\lbiprod H_F$.\endproof

Next we wish to show that the family $\widetilde \n$ has topological charge
equal to $-1$. Recall that a pair of projections $P$, $Q$ are said
to be Murray-von Neumann equivalent if there exists a partial
isometry $V$ such that $P=VV^*$ and $Q=V^*V$.

\begin{lem}The projections $1\otimes \pp$ and $\widetilde \Pp$ are Murray-von Neumann equivalent in the algebra
$\M_4((\B(\SL_\theta(2,\HH))\lbiprod H_F)\otimes \A_r(S^4_\theta))$
and hence have the same topological charge.\end{lem}

\proof First one shows as in \cite{lprs:ncfi} that the projections
$1\otimes\qp$ and $\widetilde \Qp$ are Murray-von Neumann equivalent
in the algebra $\M_4((\B(\SL_\theta(2,\HH))\lbiprod H_F)\otimes
\A_r(S^4_\theta))$. Indeed, defining a partial isometry
$V=(V_{kl})\in \M_4((\B(\SL_\theta(2,\HH))\lbiprod H_F)\otimes
\A_r(S^4_\theta))$ by
$$V_{kl}:=\sum_\alpha\rho^{-1}(\hatA_{k\alpha}\otimes\tau_\alpha)\otimes\qp_{\alpha
l},$$ a straightforward computation shows that $V^*V=1\otimes \qp$
and $VV^*=\widetilde \Qp$. Since $1\otimes \qp$ and $1\otimes \pp$
are complementary projections, as are $\widetilde \Qp$ and
$\widetilde \Pp$, it immediately follows that $1\otimes \pp$ and
$\widetilde \Pp$ are Murray-von Neumann equivalent. One may also
show in the same way as in \cite{lprs:ncfi} that the topological
charge of the projection $\widetilde\Qp$ is equal to $1$, whence it
follows that $\widetilde\Pp$ has topological charge equal to
$-1$.\endproof

As we did for the basic instanton, we can produce other families by
conjugating $\widetilde \Pp$ with a unitary matrix. In particular,
we can take any quadruple $(u_1,u_2,u_3,u_4)$ of unitary elements in
$H_F$ with $u_1=u_2^*$, $u_3=u_4^*$ and define
$$U=\textup{diag}\,\left(u_1,u_2,u_3,u_4\right)\in\M_4((\B(\SL_\theta(2,\HH))\lbiprod
H_F)\otimes\A_r(S^4_\theta))$$ (we have suppressed the trivial
factors in the tensor product in the entries of $U$). In this case,
the conjugated projection $U\widetilde \Qp U^*$ has entries
\begin{align}\label{big uni}(U\widetilde \Qp U^*)_{kl}&=U_{k\alpha}\widetilde
\Pp_{\alpha\beta}(U_{l\beta})^*\\\nonumber
&=\sum_{\alpha,\beta}\rho^{-2}\hatA_{k\alpha}\,\underline{\cdot}\,(\hatA^*)_{\beta
l}\otimes u_k\tau_\alpha\tau_\beta^*u_l^*\otimes
(\sfu\sfu^*)_{\alpha\beta}\,F^{-2}(\tau_\beta\tau_l^*,\tau_\alpha\tau_\beta^*),\end{align}
which are elements in the algebra $(\B(\SL_\theta(2,\HH))\lbiprod
H_F)\otimes\A_r(S^4_\theta)$. The conjugated projection
$U\widetilde\Pp U^*$ yields a family of instantons parameterised by
the algebra $\B(\SL_\theta(2,\HH))\lbiprod H_F$ which is gauge
equivalent to the Grassmann family defined by $\widetilde\Pp$.

Finally in this section, we also consider the coaction of the Hopf
algebra $\B(\Sp_\theta(2))\lbiprod H_F$ on the basic instanton
defined by $\pp$, with the following result.
We denote by $\widetilde\Pp_0$ and $\widetilde \Qp_0$ the images of
the projections $\pp$ and $\qp$ under the coaction
$\A(S^4_\theta)\to \left(\B(\Sp_\theta(2))\lbiprod H_F\right)\otimes
\A(S^4_\theta)$.

\begin{prop}\label{uni coact}The Grassmann connection $\widetilde\n_0:=\widetilde
\Pp_0\circ(\id\otimes\D)$ is a family of instantons parameterised by
the Hopf algebra $\B(\Sp_\theta(2))\lbiprod H_F$; the family
$\widetilde\n_0$ is equivalent to the basic instanton
$\n:=\pp\circ\D$.\end{prop}

\proof The fact that $\widetilde\n_0$ is a family of instantons
follows in the same way as the proof of Proposition~\ref{big fam}.
We have to show that the projection $\widetilde\Pp_0$ is unitarily
equivalent to $1\otimes\pp$ in the matrix algebra
$\M_4((\B(\Sp_\theta(2))\lbiprod H_F)\otimes \A(S^4_\theta))$. The
unitary matrix which achieves this is $$U=(U_{kl})\in
\M_4((\B(\Sp_\theta(2))\lbiprod H_F)\otimes \A(S^4_\theta)),\qquad
U_{kl}:=\big(\hatA_{kl}\otimes \tau_l\big)\otimes 1,$$ where the
elements $\hatA_{kl}$ in this case denote the generators of
$\B(\Sp_\theta(2))$. It is a straightforward computation to check
that $\widetilde\Qp_0=U(1\otimes \qp)U^*$ and hence that
$\widetilde\Pp_0=
U(1\otimes\pp)U^*,$ whence the result. In these computations
it is important to note that, thanks to the $*$-structure on
$\B(\Sp_\theta(2))\lbiprod H_F$, the matrix $U^*$ has entries
$(U^*)_{kl}=F^{-2}(\tau_k\tau_l^*,\tau_k^*)(\hatA^*)_{kl}\otimes\tau_k^*\otimes
1$.\endproof

In summary, we have shown how to construct various
families of instantons on $S^4_\theta$ all having topological charge
equal to $-1$. Clearly, the parameter spaces for these families are
not necessarily `optimal', in the sense that some of the parameters
may correspond to gauge equivalent instantons. It is of central
interest and of course only natural to investigate how to remove
these extra gauge parameters and leave parameter spaces which
describe gauge equivalence classes of instantons. This is addressed
in the next section.

\section{Parameter Spaces for Charge One Instantons}\label{sect charge one params} 


In the classical case, if a Lie group $G$ acts
(freely, let us say) on a smooth manifold $P$, then one can consider the
corresponding space of orbits $P/G$. 
In cases when $P$ is a parameter space for a
family of instantons such that the action of $G$ 
preserves gauge equivalence classes, then one can always construct a new family of
instantons labeled by the more ``efficient''
parameter space $P/G$ (the latter clearly having less redundancy).

In the noncommutative setting, where we allow for noncommutative
parameter spaces, our strategy is analogous, with group actions and
spaces being replaced by coactions of Hopf algebras and
`noncommutative quotients'. We show that, in the situation
where the coaction of a Hopf algebra on a parameter space results in
a gauge equivalent family of instantons, there is a family of
instantons parameterised by the noncommutative quotient space (the
algebra of coinvariants for the coaction).

As mentioned in Remark~\ref{def mod}, we think of the moduli functor
as a functor whose source is the category ${}^{H_F}\textsf{Alg}$ of
unital left $H_F$-comodule $*$-algebras, {\em i.e.} we consider
parameter spaces described by algebras which carry a left
$H_F$-coaction. This is in keeping with our strategy of viewing the
passage from classical to quantum as a `quantisation functor' as in
\S\ref{section cocycle twists}. In particular, this means that all
quantum groups we consider are braided Hopf algebras in the
category, and all coactions are required to be morphisms in the
category and hence braided as well.

\subsection{Removing the $\B(\Sp_\theta(2))$ gauge parameters}\label{params}
For the sake of brevity, in this section we write
$A:=\B(\SL_\theta(2,\HH))\lbiprod H_F$ and $C:=\B(\Sp_{\theta}(2))$.
Of all the families of charge one instantons that we
have constructed, the largest is the one parameterised by the
noncommutative algebra $A$, and we seek to make it smaller by
quotienting away the parameters corresponding to gauge equivalence.
In this section, we consider a coaction of $C$ on the parameter
space $A$ and show that it generates a gauge equivalent family of
instantons. These gauge parameters are removed by constructing the
corresponding quantum quotient of $A$ by $C$.

Recall that the braided Hopf algebra $\B(\Sp_\theta(2))$ is the
quotient of $\B(\SL_\theta(2,\HH))$ by an ideal
$\mathcal{I}_\theta$, obtained as a twist of the ideal $\mathcal{I}$
defined in Eq.~\eqref{eqn ideal gens}. Let us write
$$\Pi_\theta:\B(\SL_\theta(2,\HH))\rightarrow \B(\Sp_\theta(2))$$
for the canonical projection. Using this we can define a braided
left coaction of $C$ on $A$ as follows. Note that $A$ is canonically
an object in the category ${}^{H_F}\mM$ {\em via} the tensor product
$H_F$-coaction, whence we may form the braided tensor product
algebra $C\utimes A$. We use the notation $\Delta_F(a)=a\uo\otimes
a\ut$ for the braided coproduct $\Delta_F$ of
$\B(\SL_\theta(2,\HH))$.

\begin{lem}There is a braided left coaction of the braided Hopf algebra $C$ on the parameter space $A$ defined by the formula $$\delta_L:A\to C\utimes A,\qquad\delta_L(a\otimes h)=\Pi_\theta(a\uo)\otimes a\ut\otimes h,$$
for which $A$ is a braided left $C$-comodule algebra.\end{lem}

\proof The fact that $\delta_L$ defines a braided coaction is
immediate from the fact that both the coproduct $\Delta_F$ and the
projection $\Pi_\theta$ are morphisms in the category ${}^{H_F}\mM$,
hence so is the composition $\delta_L=\Pi_\theta \circ \Delta_F$. To
check that $\delta_L$ respects the algebra structure of $A$, we
compute on generators that
\begin{align*}\delta_L(\hatA_{ij}\otimes
h)\,\underline{\cdot}\,\delta_L(\hatA_{kl}\otimes
g)&=\sum_{\alpha,\beta}(\Pi_\theta(\hatA_{i\alpha})\otimes\hatA_{\alpha
j}\otimes
h)\,\underline{\cdot}\,(\Pi_\theta(\hatA_{k\beta})\otimes\hatA_{\beta
l}\otimes
g)\\&=\sum_{\alpha,\beta}\Pi_\theta(\hatA_{i\alpha}\,\underline{\cdot}\,\hatA_{k\beta})\otimes\hatA_{\alpha
j}\,\underline{\cdot}\,\hatA_{\beta l}\otimes hg \times
\\&\qquad\qquad\qquad\qquad\qquad\times F^{-2}(\tau_k\tau_\beta^*,\tau_\alpha\tau_j^*h)F^{-2}(\tau_\beta\tau_l^*,h)\\&=\sum_{\alpha,\beta}\Pi_\theta(\hatA_{i\alpha}\,\underline{\cdot}\,\hatA_{k\beta})\otimes\hatA_{\alpha
j}\,\underline{\cdot}\,\hatA_{\beta l}\otimes hg \times
\\&\qquad\qquad\qquad\qquad\qquad\times
F^{-2}(\tau_k\tau_l^*,h)F^{-2}(\tau_k\tau_\beta^*,\tau_\alpha\tau_j^*)\\&=\delta_L(\hatA_{ij}\,\underline{\cdot}\,\hatA_{kl}\otimes
hg)F^{-2}(\tau_k\tau_l^*,h)\\&=\delta_L((\hatA_{ij}\otimes
h)\,\underline{\cdot}\,(\hatA_{kl}\otimes g)),\end{align*} as
required for a braided comodule algebra.\endproof

This establishes the quantum analogue of a group action on our
parameter space. The following lemma tells us that this action is by
gauge transformations. Since our strategy is to compare pairs of
parameter spaces by looking to see if they define unitarily
equivalent families of instantons (through an application of the
moduli functor defined in Remark~\ref{def mod}), the correct way to
interpret the effect of the coaction $\delta_L$ on the family
defined by the parameter space $A$ is by coacting upon the entries
of the projection $\Pp \in\M_4\big(A\otimes\A_r(S^4_\theta)\big)$ by
$\delta_L\otimes\id$, that is to say by leaving the algebra
$\A_r(S^4_\theta)$ alone in this coaction.

\begin{prop}
The image $\delta_L(\widetilde\Pp)$ of the projection
$\widetilde\Pp$ under the coaction $\delta_L\otimes\id$ is unitarily
equivalent to the projection $1\otimes\widetilde\Pp$ in the algebra
$\M_4(\big(C\utimes A\big)\otimes\A_r(S^4_\theta))$.
\end{prop}

\proof We first consider the effect of the coaction
$\delta_L\otimes\id$ on the projection $\widetilde \Qp$ in \eqref{qtproj}:
\begin{align*}
(\delta_L\otimes\id)
(\widetilde \Qp_{kl})&=\sum_{\alpha,\beta,\gamma,\delta}(1\otimes\rho^{-2})\Pi_\theta(\hatA_{k\gamma}\,\underline{\cdot}\,(\hatA^*)_{\delta
l})\otimes\hatA_{\gamma\alpha}\,\underline{\cdot}\,(\hatA^*)_{\beta\delta}\otimes\tau_\alpha\tau_\beta^*\otimes
(\sfu\sfu^*)_{\alpha\beta}\,\times
\\&\qquad\qquad\qquad\times F^{-2}(\tau_\beta\tau_l^*,\tau_\beta^*)
F^{-2}(\tau_\delta\tau_l^*,\tau_\beta\tau_\delta^*)F^{-2}(\tau_\delta\tau_l^*,\tau_\gamma)F^{-2}(\tau_\beta\tau_\delta^*,\tau_\alpha).\end{align*}It
is a straightforward calculation, along the same lines as
Proposition~\ref{uni coact}, to check that the same effect is
achieved by conjugation with the unitary matrix
$$U=(U_{kl})\in\M_4(\big(C\utimes A\big)\otimes\A_r(S^4_\theta)),\qquad
U_{kl}=\big(\Pi_\theta(\hatA_{kl})\otimes 1\big)\otimes 1,$$ {\em
i.e.} we have $\delta_L(\widetilde \Qp)=U(1\otimes\widetilde
\Qp)U^*$. From the fact that $\widetilde\Pp$ is complementary to
$\widetilde \Qp$ it follows immediately that
$\delta_L(\widetilde\Pp)=U(1\otimes\widetilde \Pp)U^*$, as
required.\endproof

The subalgebra $C=\B(\Sp_\theta(2))$ of
$A=\B(\SL_\theta(2,\HH))\lbiprod H_F$ thus consists entirely of
gauge parameters which we would like to remove. This reduction of
parameters is performed by computing the quantum quotient of $A$ by
the coaction of $C$.

\begin{prop}The algebra of coinvariants for the coaction $\delta_L$,
$$\{a \in
\B(\SL_\theta(2,\HH))\lbiprod H_F~|~\delta_L(a)=1\otimes a\},
$$
 is isomorphic to the subalgebra
$\B(\mathsf{M}_{\theta})\lcross H_F$ of
$\B(\SL_\theta(2,\HH))\lbiprod H_F$, where $\B(\mathsf{M}_{\theta})$
is the subalgebra of $\B(\SL_\theta(2,\HH))$ generated by the
elements $$\mathsf{m}_{ij}:=\sum_\alpha(\hatA_{\alpha
i})^*\,\underline{\cdot}\,\hatA_{\alpha j}, \qquad i,j=1,\ldots,4.$$
\end{prop}

\proof Since the relations in $\B(\Sp_\theta(2))$ are quadratic in
the generators $\hatA_{kl}$ and their conjugates, the generators of
the algebra of coinvariants must be at least quadratic in them. The
key relations here are those coming from the antipode, namely
$$\sum_\alpha(A^*)_{i\alpha}A_{\alpha j}=\delta_{ij}, \qquad
i,j=1,\ldots,4.$$ In order for the first leg of the tensor product
in $\delta_L(a\otimes h)$ to involve these relations, we have to
take $a=\sum_\alpha(\hatA_{\alpha
i})^*\,\underline{\cdot}\,\hatA_{\alpha j}$. Moreover, we compute
that for all group-like elements $h\in H_F$ we have
\begin{align*} \delta_L\left(\sum_\alpha(\hatA_{\alpha i})^*\,\underline{\cdot}\,\hatA_{\alpha j}\otimes h\right)&=\sum_{\alpha,\beta,\gamma}\Pi_\theta\left((\hatA_{\alpha \beta})^*\,\underline{\cdot}\,\hatA_{\alpha \gamma}\right)
\otimes (\hatA_{\beta i})^*\,\underline{\cdot}\,\hatA_{\gamma
j}\otimes
h\\&=\sum_{\beta,\gamma}\delta_{\beta\gamma}\otimes(\hatA_{\beta
i})^*\,\underline{\cdot}\,\hatA_{\gamma j} \otimes h
\\&=1\otimes \sum_\beta((\hatA_{\beta i})^*\,\underline{\cdot}\,\hatA_{\beta j})\otimes h.\end{align*} The identification of the algebra of coinvariants as $\B(\mathsf{M}_\theta)\lcross H_F$ is now obvious.\endproof

This gives us a new parameter space which we denote by
$B:=\B(\mathsf{M}_{\theta})\lcross H_F$. Next we have to check that
it really does parameterise a family of instantons. To this end, let
$\E_A$ denote the finitely-generated projective
$A\otimes\A_r(S^4_\theta)$-module defined by the projection
$\widetilde \Pp$. We define
$$\E_B=\E_A^{\textup{co}\,C}:=\{\xi\in\E_A~|~\delta_L(\xi)=1\otimes \xi\}$$ to be the vector space of coinvariant elements in
$\E_A$. It is clear that the right $A\otimes\A_r(S^4_\theta)$-module
structure on $\E_A$ survives as a right
$B\otimes\A_r(S^4_\theta)$-module structure on $\E_B$.

\begin{lem}\label{can BA iso}The induced module $$\E_B\otimes_{B\otimes\A_r(S^4_\theta)}(A\otimes \A_r(S^4_\theta))\simeq
\E_B\otimes_BA$$ is canonically isomorphic to $\E_A$ as a right
$A\otimes\A_r(S^4_\theta)$-module.\end{lem}

\proof The proof of this assertion goes along the lines of
\cite{schn}, noting that the proof there is given in terms of
`ordinary' rather than braided coactions. However, the proof goes
through in the braided case as well: in fact one may view the
inclusion $B\hookrightarrow A$ purely as an extension of coalgebras
and still make the same conclusions \cite{tb}, so that the (braided
or ordinary) algebra structure of the extension is not important.
The strategy is to check that the canonical algebra inclusion
$\iota:B\hookrightarrow A$ is a faithfully flat (braided)
$C$-Hopf-Galois extension: this follows from the facts that the
canonical linear map
\begin{equation}\label{can map}\chi:A\otimes_BA\to C\otimes A,
\qquad a\otimes a'\mapsto \delta_L(a)a'\end{equation} is a bijection
(as is always the case for coactions defined in this way by a Hopf
algebra projection) and that $C=\B(\Sp_\theta(2))$ is a cosemisimple
Hopf algebra. From this, it follows that the category
${}^C\mM_{A\otimes\A_r(S^4_\theta)}$ of left $C$-comodule right
$A\otimes\A_r(S^4_\theta)$-modules is equivalent to the category
$\mM_{B\otimes\A_r(S^4_\theta)}$ of right
$B\otimes\A_r(S^4_\theta)$-modules (each viewed as a sub-category of
${}^{H_F}\mM$).\endproof

The module $\E_B$ thus defines a family of Hermitian vector bundles
over $S^4_\theta$ parameterised by the algebra $B$. On the
projective $A\otimes \A_r(S^4_\theta)$-module $\E_A=\widetilde
\Pp(A\otimes\A_r(S^4_\theta))^4$ we have the family of instantons,
$\n_A:=\widetilde \Pp\circ(\id\otimes \D)$, as constructed in
\S\ref{nc family}. The next proposition gives us the required family
of instanton connections on the family of bundles $\E_B$.

\begin{prop}\label{quotient}Let $\iota:B\hookrightarrow A$ be the canonical algebra inclusion. Then there exists a Grassmann family $\n_B$ of
instantons parameterised by the algebra $B$, unique up to unitary
equivalence, with the property that $\iota_*(\n_B)=\n_A$.\end{prop}

\proof Recall that $\n_A$ is the Grassmann connection on the
projective $A\otimes \A_r(S^4_\theta)$-module $\E_A=\widetilde
\Pp(A\otimes\A_r(S^4_\theta))^4$. From the above discussion, the
coinvariant sub-module $\E_B:=\E_A^{\textup{co}\,C}$ is
finitely-generated and projective as a right
$B\otimes\A_r(S^4_\theta)$-module, and hence defined by a projection
$\Pp_B$, unique up to unitary equivalence. We define a Grassmann
family of connections on $\E_B$ by $\n_B:=\Pp_B\circ (\id\otimes
\D)$. Since the induced module $\E_B\otimes_BA$ is canonically
isomorphic to $\E_A$, the induced family of connections
$\iota_*(\n_B)$ defined in Proposition~\ref{mod func} must be the
same as the family $\n_A$ (up to equivalence). From the proof of
Proposition~\ref{inj curv}, we see that the curvature of $\n_B$ is
the same as the curvature of $\n_A$, which means that the curvature
of the family $\n_B$ must also be anti-self-dual.\endproof

\subsection{Removing the $H_F$ gauge parameters}\label{H params}
We have thus removed the gauge parameters corresponding to the
braided group $\B(\Sp_\theta(2))$, yielding a family of instantons
parameterised by the algebra $B$. The next step is to quotient away
the parameters corresponding to the algebra $H_F$.

Our strategy is the same as before: we remove these parameters by
considering a (braided) coaction of $H_F$ on $B$. By showing that
this coaction is by unitary gauge transformations, we then pass to
the parameter space described by the quantum quotient of $B$ by
$H_F$. This a very delicate process, however, since there are many
different ways in which $H_F$ can coact upon the parameter space
$B=\B(\mathsf{M}_\theta)\lcross H_F$, whence there are many
different ways in which we can make the quotient.

Once again we consider all of our parameter spaces as being objects
in the category ${}^{H_F}\textsf{Alg}$, in particular noting that
$H_F$ is canonically an object in the category {\em via} the left
regular coaction. This means that we can form the braided tensor
product algebras $H_F\utimes A$ and $H_F\utimes B$ as objects in the
category.

\begin{lem}\label{alt torus act}Let $u:=(u_1,u_2,u_3,u_4)$ be unitary elements of $H_F$
such that $u_1^*=u_2$, $u_3^*=u_4$. Then there is a left coaction
$\delta_u:A\to H_F\utimes A$ defined on by
$$\delta_u(\hatA_{kl}\otimes h)=u_ku_l^*h\otimes \hatA_{kl}\otimes
h,\qquad k,l=1,\ldots,4$$ for each group-like element $h\in H_F$ and
extended as a braided $*$-algebra map.\end{lem}

\proof We check the conditions for $\delta_u$ to define a braided
coaction of $H_F$:
\begin{align*}\big((\id\otimes\delta_u)\circ\delta_u\big)(\hatA_{kl}\otimes
h)&=(\id\otimes\delta_u)(u_ku_l^*h\otimes \hatA_{kl}\otimes
h)\\&=u_ku_l^*h\otimes u_ku_l^*h\otimes \hatA_{kl}\otimes
h\\&=\big((\Delta\otimes \id)\circ\delta_u\big)(\hatA_{kl}\otimes
h);\\\big((\ep\otimes\id)\circ\delta_u\big)(\hatA_{kl}\otimes
h)&=\ep(u_ku_l^*h)\hatA_{kl}\otimes h=\hatA_{kl}\otimes
h,\end{align*} where in the last line we have used the fact that
$\ep(h)=1$ for all group-like elements $h\in H_F$. We then extend
$\delta_u$ as a braided $*$-algebra map to obtain the
result.\endproof

As before, we extend the coaction $\delta_u:A\to H_F\utimes A$ to a
coaction on the algebra $A\otimes\A(S^4_\theta)$ by
$\delta_u\otimes\id$. In this way we can coact upon the projection
$\widetilde\Pp$ with $H_F$. Our next result is that this coaction is
by gauge transformations.

\begin{prop}\label{del equiv}Let $u:=(u_1,u_2,u_3,u_4)$ be unitary elements of $H_F$ as above. Then the image $(\delta_u\otimes\id)(\widetilde\Pp)$ of the projection
$\widetilde\Pp$ under the coaction $\delta_u\otimes\id$ is unitarily
equivalent to the projection $1\otimes\widetilde\Pp$ in the algebra
$\M_4(\big(H_F\utimes A\big)\otimes\A_r(S^4_\theta))$.\end{prop}

\proof We first consider the effect of the coaction $\delta_u$ on
the projection $\widetilde \Qp$:
\begin{align*}(\delta_u\otimes\id)(\widetilde
\Qp_{kl})&=(1\otimes\rho^{-2})\sum_{\alpha,\beta}
\big(u_ku_\alpha^*u_\beta
u_l^*\tau_\alpha\tau_\beta^*\big)\otimes\hatA_{k\alpha}\,\underline{\cdot}\,(\hatA^*)_{\beta
l}\otimes\tau_\alpha\tau_\beta^*\otimes (\sfu\sfu^*)_{\alpha\beta}\,
\times
\\&\qquad\qquad\qquad\qquad\qquad\qquad \times F^{-2}(\tau_\beta\tau_l^*,\tau_\alpha\tau_\beta^*).\end{align*}
It is a straightforward calculation to check that the same effect
can be achieved by conjugating with the unitary diagonal matrix
$$U=(U_{kl})\in\M_4(\big(H_F\utimes A\big)\otimes\A_r(S^4_\theta)),\qquad U_{kl}=\textup{diag}\,\big((u_k\tau_k\otimes 1)\otimes 1~|~k=1,\ldots,4\big),$$ {\em
i.e.} we have $(\delta_u\otimes\id)(\widetilde
\Qp)=U(1\otimes\widetilde \Qp)U^*$. From the fact that
$\widetilde\Pp$ is complementary to $\widetilde \Qp$ it follows
immediately that
$(\delta_u\otimes\id)(\widetilde\Pp)=U(1\otimes\widetilde \Pp)U^*$,
as required.\endproof

Each of the coactions $\delta_u$ therefore gives us an equally valid
way of quotienting the parameter space and removing gauge freedom.
It is clear that the coaction $\delta_u$ descends to a coaction on
the algebra $B=\B(\mathsf{M}_\theta)\lcross H_F$, whence we have the
following result.

\begin{prop}Let $u:=(u_1,u_2,u_3,u_4)$ be unitary elements of $H_F$ as above. Then the algebra of coinvariants
$$\A(\mathsf{M}^u_\theta):=\{b\in B~|~\delta_u(b)=1\otimes b\}$$
for the coaction $\delta_u$ is generated by the elements
$M_{ij}^u:=\mathsf{m}_{ij}\otimes u_i^*u_j$ for
$i,j=1,\ldots,4$.\end{prop}

\proof This is a matter of noticing that on generators
$\mathsf{m}_{ij}\otimes h$ of $B$ the coaction $\delta_u$ has the
form
$$\delta_u(\mathsf{m}_{ij}\otimes
h)=u_iu_j^*h\otimes\mathsf{m}_{ij}\otimes h.$$ This implies that the
algebra of coinvariants is generated by elements for which
$h=u_i^*u_j$, as claimed.\endproof

By removing gauge parameters, we have thus produced a more efficient
parameter space $M:=\A(\mathsf{M}^u_\theta)$. Again we have to check
that $M$ really does parameterise a family of instantons. This time
we define
$$\E_M=\E_B^{\textup{co}\,H_F}:=\{\xi\in\E_B~|~\delta_u(\xi)=1\otimes \xi\}$$ to be the vector space of coinvariant elements in
$\E_B$. In this case the right $B\otimes\A(S^4_\theta)$-module
structure on $\E_A$ survives as a right
$M\otimes\A_r(S^4_\theta)$-module structure on $\E_M$.

\begin{lem}The induced module $$\E_M\otimes_{M\otimes\A_r(S^4_\theta)}(B\otimes \A_r(S^4_\theta))\simeq
\E_M\otimes_MA$$ is canonically isomorphic to $\E_B$ as a right
$B\otimes\A_r(S^4_\theta)$-module.\end{lem}

\proof The strategy is the same as in Lemma~\ref{can BA iso}, in
that the result follows from showing that the canonical algebra
inclusion $\iota:M\hookrightarrow B$ is a faithfully flat braided
Hopf-Galois extension. Associated to the coaction $\delta_u$ we have
the corresponding canonical linear map
$$\chi_u:B\otimes_{M}B\to H_F\otimes B,\qquad
b\otimes b'\to\delta_u(b)b',$$ which we would like to show is a
bijection. Since $H_F$ is cosemisimple as a coalgebra, it is
sufficient to check that $\chi_u$ is surjective \cite{schn}.
Moreover, it is known that $\chi_u$ is surjective if, whenever $h$
is a generator of $H_F$, then the element $1\otimes h$ is in its
image \cite{ps:gal}. The canonical map here works out on
finitely-generated elements of $H_F$ to be
\begin{align*}\chi_u\left((1\otimes h)\otimes(1\otimes h')\right)&=1\otimes hh'\otimes
(h')^*,\end{align*}
so that in order to find $\tau_j$ in the image of
$\chi_u$ for each $j=1,\ldots,4$ we can simply take $h=\tau_j$ and
$h'=\tau_j^*$.\endproof

The module $\E_M$ defines a family of Hermitian vector bundles over
$S^4_\theta$ parameterised by the algebra
$M=\A(\mathsf{M}_\theta^u)$. As before, there is a corresponding
family of instantons parameterised by this space.

\begin{prop}\label{H gauge}Let $\iota:M\hookrightarrow B$ be the canonical algebra inclusion. There exists a Grassmann family of instantons $\n_u$ parameterised by the algebra $\A({\sf M}^u_\theta)$, unique up to unitary equivalence, with the property that $\iota_*(\n_u)=\n_B$.\end{prop}

\proof This follows in exactly the same way as the proof of
Proposition~\ref{quotient}. The coinvariant submodule $\E_M$ is
finitely-generated and projective, hence defined by a projection
$\Pp_u$. We define the required family $\n_u$ of instantons by
$\n_u:=\Pp_u\circ(\id\otimes\D)$.\endproof

\begin{prop}Let $u=(u_1,u_2,u_3,u_4)$ and $v=(v_1,v_2,v_3,v_4)$ be quadruples of
unitary elements in $H_F$ such that $u_1^*=u_2$, $u_3^*=u_4$,
$v_1^*=v_2$, $v_3^*=v_4$. Then the families of instantons $\n_u$ and
$\n_v$ described by the parameter spaces $\A({\sf M}^u_\theta)$ and
$\A({\sf M}^v_\theta)$ are gauge equivalent.\end{prop}

\proof The families $\n_u$ and $\n_v$ are defined as in
Proposition~\ref{H gauge}, with corresponding parameter spaces
$\A({\sf M}^u_\theta)$ and $\A({\sf M}^v_\theta)$ arising as
coinvariant subalgebras for the coactions $\delta_u$ and $\delta_v$
of $H_F$ on $B=\B({\sf M}_\theta) \lcross H_F$. These parameter
spaces each sit inside $B$ {\em via} the canonical algebra
inclusions $\iota_u:\A({\sf M}^u_\theta)\hookrightarrow B$ and
$\iota_v:\A({\sf M}^v_\theta)\hookrightarrow B$. The result follows
from the fact that the coactions $\delta_u$ and $\delta_v$ are
themselves unitarily equivalent, as one may infer from
Eq.~\eqref{big uni}, for example.\endproof

\subsection{The space ${\sf M}^u_\theta$ of connections}We would like to compute the algebra $\A({\sf M}^u_\theta)$
explicitly. For convenience, we restrict our attention to the
following special case. Let $(r_1,r_2)\in \ZZ^2$ be a pair of
integers and take
$$u=(u_1,u_2,u_3,u_4):=(\tau_1^{r_1},\tau_2^{r_1},\tau_3^{r_2},\tau_4^{r_2}).$$ Then for this $u$, we can arrange the generators $M^u_{ij}$ of the
algebra $\A(\sfM^u_\theta)$ into a matrix $M^u_\theta=(M^u_{ij})$.
Explicitly, one finds that
$$M^u_\theta=\begin{pmatrix}m&0&g_1&g_2^*\\0&m&-\bar\nu g_2&\nu g_1^*\\g_1^*&-\nu g_2^*&n&0\\g_2&\bar\nu g_1&0&n\end{pmatrix},\qquad \nu:=\mu^{r_1-r_2+1},$$
where $m:=M^u_{11}$, $n:=M^u_{33}$, $g_1:=M^u_{13}$, $g_2:=M^u_{41}$
and $\mu:=e^{\ii \pi\theta}$ is the deformation parameter. The
relations between these generators depend of course on the choice of
integers $r_1$, $r_2$. We compute them as follows.

\begin{prop}The relations between the entries of the matrix $M^u_\theta$ in the algebra
$\A({\sf M}^u_\theta)$ are given by
$$g_1g_2=\nu^2 g_2g_1,\quad g_1g_2^*=\bar\nu^2g_2^*g_1, \quad g_1g_1^*=g_1^*g_1,\quad g_2g_2^*=g_2^*g_2$$ and $m$, $n$
central. There is also a quadric relation $$mn-\bar\nu\mu
g_1^*g_1+\nu\bar\mu g_2^*g_2=1.$$\end{prop}

\proof Computing the commutation relations is a simple calculation
using the relations in the algebra $A=\B(\SL_\theta(2,\HH))\lbiprod
H_F$. For the quadric relation, one computes that
$mn=\sfm_{11}\sfm_{33}\otimes 1$,
$g_1^*g_1=\mu^{r_2-r_1}\sfm_{31}\sfm_{13}\otimes 1$ and
$g_2^*g_2=\mu^{r_1-r_2}\sfm_{41}\sfm_{14}\otimes 1$, whence we have
that $$mn-\bar\nu\mu g_1^*g_1+\nu\bar\mu
g_2^*g_2=\left(\sfm_{11}\sfm_{33}-\sfm_{31}\sfm_{13}+\sfm_{41}\sfm_{14}\right)\otimes
1=\det(A_\theta)\otimes 1,$$ where $A_\theta$ is the
$\theta$-deformed version of the matrix in Eq.~\eqref{eqn defining
M(H)}. The relation
$$\det(A_\theta)=\sfm_{11}\sfm_{33}-\sfm_{31}\sfm_{13}+\sfm_{41}\sfm_{14}$$
is computed as in \cite{lprs:ncfi}. The fact that $\det(A_\theta)=1$
in $\B(\SL_\theta(2,\HH))$ gives the relation as stated.\endproof

We see that the choice of unitary $u$ affects both the commutation
relations in the algebra $\A({\sf M}^u_\theta)$ as well as the
quadric relation. We emphasise the following two important cases.

\begin{example}Any choice for which $r_1=r_2$ recovers the noncommutative parameter space
discovered in \cite{lprs:ncfi}, whose algebra relations have $m,n$
central and $g_1g_2=\mu^2 g_2g_1$, $g_1g_2^*=\bar\mu^2 g_2^*g_1$.
The quadric relation is the same as the classical one, $mn-
g_1^*g_1+g_2^*g_2=1$.\end{example}

\begin{example}For any choice which has
$r_2=r_1+1$ we have $\nu=1$ and hence we obtain a commutative
parameter space, {\em i.e.} the generators $m,n,g_1,g_2$ and their
conjugates all commute. However, the quadric relation in this case
is deformed, $mn-\mu g_1^*g_1+\bar\mu g_2^*g_2=1$.\end{example}

It is not difficult to see that, in the classical limit, the
algebras $\A({\sf M}^u_\theta)$ describe different but gauge
equivalent parameterisations of the same space. In the
noncommutative case, these parameter spaces are evidently different,
some being noncommutative and others classical, but they are
nevertheless still all gauge equivalent.

\section{Instantons with Higher Topological Charge}\label{sect adhm}
In this section we generalise the previous construction to treat parameter
spaces for instantons of higher topological charge. In
\cite{bl:adhm} we gave a deformed version of the ADHM construction
which produced noncommutative families of instantons with arbitrary
charge. We start with a review of this construction which emphasises
how it may be viewed in the context of braided geometry. As we did
for the charge one case, we then show how to use gauge theory to
recover commutative parameter spaces.

\subsection{A noncommutative space of monads}\label{nc monads}We begin with a description of the space of monads over the
classical space $\C^4$, which we shall later deform by means of the
twisting cocycle $F$ on the torus algebra $H=\A(\TT^2)$ that has been used
throughout this paper. We adopt the categorical approach
used in the charge one case and demand that all of our constructions
are $H$-covariant; the quantisation functor will produce the twisted version.

The algebra $\A(\C^4)$ has a natural $\ZZ$-grading given on
generators by
$$\textup{deg}(z_j)=1, \quad \textup{deg}(z_j^*)=-1, \qquad
j=1,\ldots,4.$$  This gives rise to a decomposition into homogeneous
subspaces $\A(\C^4)=\bigoplus_{n\in \ZZ}\A_n$. For each $r \in \ZZ$
we denote by $\A(\C^4)(r)$ the `degree shifted' algebra, whose
degree $n$ component is defined to be $\A_{n+r}$. Similarly, for
each finite dimensional vector space $\mH$ the corresponding free
right module $\mH \otimes \A(\C^4)$ is $\ZZ$-graded by the grading
on $\A(\C^4)$, and the shift maps on $\A(\C^4)$ induce shift maps on
$\mH \otimes \A(\C^4)$.

\begin{defn}\label{monad}Let $k\in\ZZ$ be a fixed positive integer. A {\em monad} over
the algebra $\A(\C^4)$ is a sequence of free right
$\A(\C^4)$-modules,
\begin{equation} \label{eqn module monad}0\to \mH\otimes \A(\C^4)(-1) \xrightarrow{\sigma_z}
\mK\otimes \A(\C^4) \xrightarrow{\tau_z} \mL\otimes \A(\C^4)(1)\to
0,\end{equation} where $\mH$, $\mK$ and $\mL$ are complex vector
spaces of dimensions $k$, $2k+2$ and $k$ respectively, such that the
maps $\sigma_z$, $\tau_z$ are linear in the generators
$z_1,\ldots,z_4$ of $\A(\C^4)$. The first and last terms of the
sequence are required to be exact, so that the only non-trivial
cohomology is in the middle term.\end{defn}

As in \cite{bl:adhm}, our strategy is to find the space of all
possible monads for a fixed choice of positive integer $k$. We begin
by considering the module map
$\sigma_z
$ in the complex \eqref{eqn module monad}. Choosing ordered bases
$(u_1,\ldots,u_k)$ for the vector space $\mH$ and
$(v_1,\ldots,v_{2k+2})$ for the vector space $\mK$, we can express
$\sigma_z$ as
\begin{equation}\label{module map}
\sigma_z: u_b\otimes Z \mapsto \sum\nolimits_{a,j}
M_{ab}^j \otimes v_a \otimes z_j Z, \qquad Z\in
\A(\C^4),\end{equation} for $(2k+2)\times k$ matrices
$M^j:=(M^j_{ab})$, where $j=1,\ldots,4$ and $a=1,\ldots,2k+2$,
$b=1,\ldots,k$. Thus, in more
compact notation, $\sigma_z$ may be written
\begin{equation} \label{classical sigma}
\sigma_z=\sum\nolimits_j M^j \otimes z_j.
\end{equation}
 In dual terms, we think of the $M^j_{ab}$ as coordinate
functions on the space $\M(\mH,\mK)$ of all such maps $\sigma_z$,
with (commutative) coordinate algebra $\A(\M(\mH,\mK))$ generated by
the functions $M^j_{ab}$ for $j=1,\ldots,4$ and $a=1,\ldots,2k+2$,
$b=1,\ldots,k$. It comes equipped with the homomorphism
\eqref{module map} of right $\A(\C^4)$-modules. In this way, the
space $\M(\mH,\mK)$ in fact has the structure of an algebraic
variety: it is the spectrum of the algebra $\A(\M(\mH,\mK))$.

As mentioned above, we wish to view the construction as taking place
in the category ${}^H\mM$. The free $\A(\C^4)$-modules appearing in
the complex \eqref{eqn module monad} are automatically objects in
${}^H\mM$; we need that the maps $\sigma_z$, $\tau_z$ are morphisms.

\begin{lem}The map $\sigma_z:=\sum\nolimits_\alpha M^\alpha \otimes
z_\alpha$ is a morphism in the category ${}^{H}\mM$ if and only if
the coordinate functions $M^j_{ab}$ carry the left $H$-coaction
given on generators by
$$
M^j_{ab}\mapsto \tau_j^*\otimes M^j_{ab} ,
$$
for each $j=1,\ldots,4$ and $a=1,\ldots,2k+2$, $b=1,\ldots,k$, making the
vector space spanned by the $M^j_{ab}$ into a left $H$-comodule.
\end{lem}

\proof Upon inspection of Eq.~\eqref{module map} we see that
$\sigma_z$ cannot possibly be an intertwiner for the $H$-coactions
on $\mH\otimes \A(\C^4)$ and on $\mK\otimes\A(\C^4)$ unless we also
allow for a coaction of $H$ on the algebra $\A(\M(\mH,\mK))$ as
well. It is clear that, for $\sigma_z$ to be $H$-covariant, this
coaction needs to be as stated in the lemma. \endproof

It follows that the algebra $\A(\M(\mH,\mK))$ is an algebra in the
category ${}^H\mM$. It possesses a certain universality property
which we discuss in Appendix~\ref{families}, reinforcing our
assertion that it is the coordinate algebra of the space of all
module maps $\sigma_z$.

We may carry out the same analysis for the map $\tau_z$. We choose a
basis $(w_1,\ldots,w_k)$ for the vector space $\mL$ and consider the
map
\begin{equation}\tau_z:v_a\otimes Z\mapsto
\sum\nolimits_{b,j} N_{ba}^j \otimes w_b\otimes z_j
Z.\end{equation} Then the commutative algebra $\A(\M(\mK,\mL))$
generated by the matrix elements
$$\{N_{ba}^j~|~a=1,\ldots,2k+2,~b=1,\ldots,k,~j=1,\ldots,4\},$$ when equipped with the morphism of right $\A(\C^4)$-modules
$$\tau_z: \mK\otimes \A(\C^4) \to
\A(\M(\mK,\mL))\otimes \mL \otimes \A(\C^4)(1),$$ is the coordinate
algebra of the space of all maps $\mK\otimes \A(\C^4) \to \mL
\otimes \A(\C^4)(1)$. In compact notation, the map $\tau_z$ has the
form
\begin{equation} \label{classical tau} \tau_z=\sum\nolimits_j
N^j \otimes z_j
\end{equation} upon collecting the generators into the $k\times
(2k+2)$ matrices $N^j:=(N^j_{ba})$. For covariance we need the left
$H$-coaction on $\A(\M(\mK,\mL))$ given by $$N^j_{ba}\mapsto
\tau_j^*\otimes N^j_{ba}$$ for $j=1,\ldots,4$ and $a=1,\ldots,2k+2$,
$b=1,\ldots,k$, which makes $\A(\M(\mK,\mL))$ into a left
$H$-comodule algebra, an object in the category ${}^H\mM$.

Next we need to address the requirement that \eqref{eqn module monad} be a
complex, {\em i.e.} that the composition $\vartheta_z:=\tau_z\circ
\sigma_z$ is zero. To obtain this in a coordinate algebra framework,
we note that the space of all right module maps $\mH\otimes
\A(\C^4)(-1)\to \mL\otimes \A(\C^4)(1)$ which are quadratic in the
generators $z_1,\ldots,z_4$ is encoded by the commutative algebra
$\A(\M(\mH,\mL))$ generated by matrix elements
$$\{T^{j,l}_{cd}~|~c,d=1,\ldots,k,~j,l=1,\ldots,4\},$$ together with the right module map $$\vartheta_z:\mH\otimes
\A(\C^4)(-1)\to \A(\M(\mH,\mL))\otimes \mL\otimes \A(\C^4)(1),$$
$$\vartheta_z:u_b\otimes Z \mapsto \sum\nolimits_{j,l,d}
T^{j,l}_{db}\otimes w_d \otimes z_j z_l Z$$ with
respect to our earlier choice of bases. The identification of
$\vartheta_z$ with the composition $\tau_z\circ\sigma_z$ appears in
coordinate form as a `coproduct' or a `gluing' of rectangular matrices
\cite{mm:qumat}, {\em i.e.} as an algebra map
\begin{equation}\label{triang coprod}\Delta:\A(\M(\mH,\mL))\to\A(\M(\mK,\mL))\otimes\A(\M(\mH,\mK)),$$
$$\Delta(T^{j,l}_{cd}):=\sum\nolimits_b
N^j_{cb}\otimes M^l_{bd},\qquad
j,l=1,\ldots,4,~c,d=1,\ldots,k.\end{equation}
Therefore, requiring
that the composition be zero results in the extra relations
\begin{equation}\label{adhm eqs}\sum\nolimits_b
\left(N^j_{cb}\otimes M^l_{bd}+N^l_{cb}\otimes
M^j_{bd}\right)=0 ,
\end{equation} for all $j,l=1,\ldots,4$,
$c,d=1,\ldots,k$.

\begin{defn}We denote by $\A(\tmM_k)$ the coordinate algebra of the space of all monads
\eqref{eqn module monad}. It is the quotient of the tensor product
algebra $\A(\M(\mK,\mL))\otimes\A(\M(\mH,\mK))$ by the relations
\eqref{adhm eqs}.\end{defn}

We are now ready to pass to the
noncommutative situation. Applying the `quantisation functor'
deforms our matrix coordinate algebras according to the following.

\begin{prop}\label{twisted maps}The relations in the algebras
$\A(\M(\mH,\mK))$ and $\A(\M(\mK,\mL))$ are deformed into
$$M^j_{ab}M^l_{cd}=\eta_{lj}M^l_{cd}M^j_{ab},\qquad
N^j_{ba}N^l_{dc}=\eta_{lj}N^l_{dc}N^j_{ba}$$ for all
$j,l=1,\ldots,4$ and all $a,c=1,\ldots,2k+2$, $b,d=1,\ldots,k$. \end{prop}

\proof We apply the deformation functor described in
\S\ref{section cocycle twists}. The products of generators are
deformed into $M^j_{ab}\cdot
M^l_{cd}=F(\tau_j^*,\tau_l^*)M^j_{ab}M^l_{cd}$ and $N^j_{ba}\cdot
N^l_{dc}=F(\tau_j^*,\tau_l^*)N^j_{ba}N^l_{dc}$ respectively, from
which the relations in the deformed algebras follow as stated.
\endproof
\noindent
We denote the resulting $H_F$-covariant algebras by
$\B(\M_\theta(\mH,\mK))$ and $\B(\M_\theta(\mK,\mL))$.
In turn, the `coproduct' in Eq.~\eqref{triang coprod} is deformed into
$$\Delta_F(T^{i,j}_{cd})=\sum\nolimits_b
N^i_{cb}\otimes M^j_{bd}\,F^{-1}(\tau_i^*,\tau_j^*),$$ although the
extra factor of $F^{-1}$ can be absorbed upon redefining the
generators -- as we did in Eq.~\eqref{redef} -- and we shall
henceforth assume this has been done, without changing our notation.
As was the case for the braided conformal group in \S\ref{se:bcg},
this $\Delta_F$ now extends as a homomorphism to the braided tensor
product algebra,
$$\Delta_F:\B(\M_\theta(\mH,\mL))\to \B(\M_\theta(\mK,\mL))\utimes
\B(\M_\theta(\mH,\mK)),$$ so that imposing that the composition
$\tau_z\circ\sigma_z$ is zero now results in the deformed relations
\begin{equation}\label{nc adhm
eqs}\sum\nolimits_r(N^j_{dr}M^l_{rb}+\eta_{jl}N^l_{br}M^j_{rd})=0\end{equation}
for all $j,l=1,\ldots,4$ and all $b,d=1,\ldots,k$, just as found in
\cite{bl:adhm}.

\begin{defn}\label{defn univ monad algebra}
Define $\B(\tmM_{\theta;k})$ to be the braided tensor product
algebra
$$\B(\M_\theta(\mK,\mL))\utimes\B(\M_\theta(\mH,\mK))$$ modulo the
relations \eqref{nc adhm eqs}.
\end{defn}

We stress that the relations \eqref{nc adhm eqs} are {\em not}
commutation relations between the matrix generators: they are rather
a set of quadratic relations in the algebra.

\subsection{The noncommutative ADHM construction}
The monads described in the previous section are by themselves
insufficient for the construction of bundles over $S^4_\theta$. In
the classical case, the cohomology of a monad is naturally a
finitely-generated projective right $\A(\C^4)$-module and hence a
bundle over $\C^4$. But one needs to ensure that this bundle is the
pull-back of some bundle over $S^4$, which is achieved by equipping
the monad with certain `reality structure'; in our deformed setting
this is incarnated as a $*$-structure on the algebra
$\B(\tmM_{\theta;k})$.

For this extra structure, we use the anti-linear map
$J:\A(\C^4_\theta)\to\A(\C^4_\theta)$ defined by
\begin{equation}\label{eqn J}J(z_1,z_2,z_3,z_4):=
(-z_2^*,z_1^*,-z_4^*,z_3^*)
\end{equation} and extended as an anti-algebra homomorphism. It is
clearly a morphism in the category of left $H_F$-comodules. For each
finite-dimensional complex vector space $\mH$ this immediately gives
a free left $\A(\C^4_\theta)$-module $\mH\otimes J(\A(\C^4_\theta))$
whose module structure is defined by $Z\tr(u\otimes J(W)):=u\otimes
J(WZ)$ for each $u\in \mH$, $W,Z\in\A(\C^4_\theta)$. Dual to this,
we have the free right $\A(\C^4_\theta)$-module $\mH^*\otimes
J(\A(\C^4_\theta))^*$, where $\mH^*$ is the dual vector space to
$\mH$ and
$J(\A(\C^4_\theta))^*:=\textup{Hom}_{\A(\C^4_\theta)}(J(\A(\C^4_\theta)),\A(\C^4_\theta))$.

Introducing the conjugate matrix generators $M^j_{ab}{}^*$ and
$N^l_{cd}{}^*$, we write $(M^j{}^\dag)_{ab}=M^j_{ba}{}^*$ and
$(N^l{}^\dag)_{cd}=N^l_{dc}{}^*$. All of this gives rise to a `dual
monad'
$$0\to \mL^*\otimes
J(\A(\C^4_\theta))^\star(-1)\xrightarrow{\tau_{J(z)}^*} \mK^*\otimes
J(\A(\C^4_\theta))^*\xrightarrow{\sigma_{J(z)}^*} \mH^*\otimes
J(\A(\C^4_\theta))^*(1),$$ where $\tau_{J(z)}^*$ and
$\sigma_{J(z)}^*$ are the `adjoint' maps defined by
$$\sigma_{J(z)}^*=\sum_j M^j{}^\dag \otimes
J(z_j)^*,\qquad \tau_{J(z)}^*=\sum_j N^j{}^\dag
\otimes J(z_j)^*.$$ We impose the condition that monads should
be self-conjugate with respect to this process, resulting in the
$*$-structure \begin{equation}\label{monad star}N^1=-M^2{}^\dag,
\quad N^2=M^1{}^\dag, \quad N^3=-M^4{}^\dag, \quad
N^4=M^3{}^\dag\end{equation} on the algebra $\B(\tmM_{\theta;k})$.
Note that the involution defined in \eqref{monad star} is compatible
with the $H_F$-coaction and hence with the algebra relations.
Although the algebra relations are slightly different, this
construction is otherwise described in more detail in
\cite{bl:adhm}.

\begin{defn}We write $\B({\sf M}_{\theta;k})$ for the quotient of the
algebra $\B(\tmM_{\theta;k})$ by the $*$-relations in
Eq.~\eqref{monad star}. It is the coordinate algebra of the space of
self-conjugate monads in the category ${}^{H_F}\mM$.\end{defn}

For self-conjugate monads, the important maps are therefore the
$(2k+2)\times k$ algebra-valued matrices \begin{align*} \sigma_z &=
M^1 \otimes z_1+M^2 \otimes z_2+M^3 \otimes z_3+M^4
\otimes z_4, \\
\sigma_{J(z)} &=-M^1 \otimes z_2^* + M^2 \otimes z_1^* -M^3\otimes
z_4^* + M^4 \otimes z_3^* ,
\end{align*}
which obey the monad conditions $\sigma_{J(z)}^\star \sigma_z=0$ and
$\sigma_{J(z)}^\star\sigma_{J(z)}=\sigma_z^\star\sigma_z$.
The crucial technical condition that we need for the ADHM
construction is the following.

\begin{lem}The entries of the
matrix $\rho^2:=\sigma_z^* \sigma_z = \sigma_{J(z)}^*\sigma_{J(z)}$
commute with the entries of the matrix $\sigma_z$.\end{lem}

\proof One finds that the $(\mu,\nu)$ entry of $\rho^2$ and the
$(a,b)$ entry of $\sigma_z$ are respectively
$$(\rho^2)_{\mu\nu}=\sum\nolimits_{r,j,l}(M^j{}^\dag)_{\mu
r}M^l_{r\nu}\otimes z_j^*z_l, \qquad
(\sigma_z)_{ab}=\sum\nolimits_{s} M^s_{ab}\otimes
z_s.$$ Suppressing the summation, the relations between these
elements are computed in the braided tensor product algebra $\B({\sf
M}_{\theta;k})\utimes \A(\C^4_\theta)$ as follows:
\begin{align*}\left((M^j{}^\dag)_{\mu r}M^l_{r\nu}\otimes
z_j^*z_l\right)\left(M^s_{ab}\otimes z_s\right)
&=(M^j{}^\dag)_{\mu r}M^l_{r\nu}M^s_{ab}\otimes
z_j^*z_l z_s\,
F^{-2}(\tau_s^*,\tau_j^*\tau_l)\\&=M^s_{ab}(M^j{}^\dag)_{\mu
r}M^l_{r\nu}\otimes z_s z_j^*z_l\,
F^{-2}(\tau_s^*,\tau_j^*\tau_l)(\eta_{js}\eta_{s l})^2\\&=\left(M^s_{ab}\otimes
z_s\right)\left((M^j{}^\dag)_{\mu r}M^l_{r\nu}\otimes
z_j^*z_l\right)\\& \qquad \qquad \qquad
F^{-2}(\tau_s^*,\tau_j^*\tau_l)(\eta_{js}\eta_{s l})^2F^{-2}(\tau_s,\tau_j\tau_l^*)\\&=\left(M^s_{ab}\otimes
z_s\right)\left((M^j{}^\dag)_{\mu r}M^l_{r\nu}\otimes
z_j^*z_l\right).\end{align*}
In the first and third
equalities we have used the definition of the braided tensor
product; in the second equality we have used the algebra relations
in $\B({\sf M}_{\theta;k})$ and $\A(\C^4_\theta)$.\endproof

We slightly enlarge the matrix algebra $\M_k(\C)\otimes \big(\B({\sf
M}_{\theta;k})\utimes \A(\C^4_\theta)\big)$ by adjoining an inverse
element $\rho^{-2}$ for $\rho^2$
and combine the matrices $\sigma_z$, $\sigma_{J(z)}$
into the $(2k+2)\times 2k$ matrix
\begin{equation}\label{matU}
\sfV :=\begin{pmatrix}\sigma_z & \sigma_{J(z)}\end{pmatrix},
\end{equation}
which by the definition of $\rho^2$ obeys
$$\sfV^* \sfV=\rho^2\begin{pmatrix} \mathbb{I}_{k} & 0 \\ 0 &
\mathbb{I}_k\end{pmatrix},$$ where $\mathbb{I}_k$ denotes the
$k\times k$ identity matrix. It follows as in \cite{bl:adhm} that
the quantity
\begin{equation}\label{proj-q}
\Qp:=\sfV\rho^{-2} \sfV^*=\sigma_z\rho^{-2}\sigma_z^* +
\sigma_{J(z)}\rho^{-2}\sigma_{J(z)}^*
\end{equation}
is automatically a $(2k+2)\times(2k+2)$ projection,
$\Qp^2=\Qp=\Qp^*$, with entries in the algebra $\B({\sf
M}_{\theta;k})\utimes \A_r(S^4_\theta)$. From this we construct the
complementary projection $\Pp:=\mathbb{I}_{2k+2}-\Qp$, having
entries in the same algebra.

At this point we encounter the same technical issue that we did in
the charge one case: for $\Pp$ to define an honest family of vector
bundles as in Definition~\ref{fam of bundles}, we need a projection
with entries in an algebra of the form $A\otimes\A_r(S^4_\theta)$
(where $A$ is the parameter space), whereas the projection $\Pp$ has
entries in the {\em braided} tensor product $\B({\sf
M}_{\theta;k})\utimes \A_r(S^4_\theta)$. In the charge one case we
had a $\B(\SL_\theta(2,\HH))\utimes \A_r(S^4_\theta)$-valued
projection, from which we passed to a
$\big(\B(\SL_\theta(2,\HH))\lbiprod H_F\big)\otimes
\A_r(S^4_\theta)$-valued projection by making a cobosonisation.
Despite the fact that $\B({\sf M}_{\theta;k})$ is only an algebra
and not a Hopf algebra, we can nevertheless use the same strategy to
obtain a genuine family of vector bundles.

Indeed, we shall convert $\Pp$ into a projection with entries in the algebra
$\B({\sf M}_{\theta;k})\lcross H_F$, where the cross product is the
one defined by the canonical left action of $H_F$ on
$\B({\sf M}_{\theta;k})$ defined in Eq.~\eqref{can act} for the general case. This action is given
on generators by the formula $$h\tr
M^j_{ab}=F^{-2}(\tau_j^*,h)M^j_{ab}, \qquad h\tr
(M^j_{ab})^*=F^{-2}(\tau_j,h)(M^j_{ab})^*,\qquad h\in H_F,$$ for
$j=1,\ldots,4$ and $a=1,\ldots,2k+2$, $b=1,\ldots,k$, and it
comes from the left coaction
$$\Delta_L:\A(\C^4_\theta)\to H_F\otimes \A(\C^4_\theta),\qquad z_j\mapsto \tau_j\otimes
z_j,$$ extended as a $*$-algebra map. More generally we shall denote
the coaction on an arbitrary element $Z\in\A(\C^4_\theta)$ by
$\Delta_L(Z)=Z\bo\otimes Z\bt$. The key result that we need is the
following.

\begin{lem}There is a $*$-algebra map $$\beta:\B({\sf M}_{\theta;k})\utimes
\A(\C^4_\theta)\to \left(\B({\sf M}_{\theta;k})\lcross
H_F\right)\otimes \A(\C^4_\theta)$$ defined by $\beta(M\otimes Z)
=M\otimes Z\bo\otimes Z\bt,$ for each $M\in \B({\sf M}_{\theta;k})$ and
$Z\in \A(\C^4_\theta)$.\end{lem}

\proof We simply check that on generators we have
\begin{align*}\beta(M^j_{ab}\otimes z_l) \beta(M^r_{cd}\otimes
z_s)&=(M^j_{ab}\otimes \tau_l\otimes z_l)(M^r_{cd}\otimes
\tau_s\otimes z_s)\\&=M^j_{ab}M^r_{cd}\otimes \tau_l\tau_s\otimes
z_lz_s\,F^{-2}(\tau_r^*,\tau_l)\\&=\beta(M^j_{ab}M^r_{cd}\otimes
z_lz_s)\,F^{-2}(\tau_r^*,\tau_l)\\&=\beta\left((M^j_{ab}\otimes
z_l)(M^r_{cd}\otimes z_s)\right)\end{align*} for all
$j,l,r,s=1,\ldots,4$, showing that $\beta$ is an algebra map.
Moreover,
\begin{align*}\left(\beta(M^j_{ab}\otimes
z_l)\right)^*&=\left(M^j_{ab}\otimes \tau_l\otimes
z_l\right)^*=M^j_{ab}{}^*\otimes \tau_l^*\otimes
z_l^*\,F^{-2}(\tau_j^*,\tau_l)\\&=\beta(M^j_{ab}{}^*\otimes
z_l^*)\,F^{-2}(\tau_j^*,\tau_l)=\beta\left((M^j_{ab}\otimes
z_l)^*\right),\end{align*} so that $\beta$ respects the
$*$-structure as well.\endproof

Immediately we apply $\beta$ to the projection $\Pp$ and, since it
is an algebra map, we obtain a projection $\widetilde \Pp$ with
entries in the algebra $\left(\B({\sf M}_{\theta;k})\lcross
H_F\right)\otimes \A_r(S^4_\theta)$. In the same way as it is shown
in \cite{bl:adhm}, the projection $\widetilde \Pp$ defines a family
of rank two Hermitian vector bundles over $S^4_\theta$, together
with the family $\widetilde\n:=\Pp\circ(\id\otimes\D)$ of Grassmann
connections whose curvature is anti-self-dual. The Chern classes of
$\widetilde \Pp$ are computed to be $\ch_1(\widetilde \Pp)=0$,
$\ch_2(\widetilde \Pp)=-k$, whence we get a family of charge $k$
instantons parameterised by the algebra $\B({\sf
M}_{\theta;k})\lcross H_F$.

\subsection{Removing the $H_F$ gauge parameters}We may now apply the
strategy of \S\ref{H params} in order to remove the gauge freedom
corresponding to the Hopf algebra $H_F$ from the family $\widetilde
\n$. To do this, we have to: choose a coaction of $H_F$ on the
parameter space $\B({\sf M}_{\theta;k})\lcross H_F$; check that this
coaction corresponds to gauge freedom; find the quotient space and
verify that it does indeed parameterise a family of instantons.

\begin{prop}Let $u:=(u_1,u_2,u_3,u_4)$ be unitary elements of $H_F$
such that $u_1^*=u_2$, $u_3^*=u_4$. Then there is a braided left
coaction $\delta_u:\B({\sf M}_{\theta;k})\lcross H_F\to H_F\utimes
(\B({\sf M}_{\theta;k})\lcross H_F)$ defined by
$$\delta_u(M^j_{ab}\otimes h)=u_j^*h\otimes M^j_{ab}\otimes
h,\qquad
$$ for group-like elements $h\in H_F$ and
extended as a braided $*$-algebra map. Moreover, the resulting
projection $\delta_u(\widetilde \Pp)$ is unitarily equivalent to the
projection $1\otimes \widetilde \Pp$ in the algebra
$\M_4\left(H_F\utimes (\B({\sf M}_{\theta;k})\lcross
H_F)\otimes\A(S^4_\theta)\right)$.\end{prop}

\proof One verifies the conditions required for $\delta_u$ to define
a braided $H_F$-comodule algebra, in the same way as was done in
Proposition~\ref{alt torus act}. The unitary equivalence is checked
in the same way as in the proof of Proposition~\ref{del
equiv}.\endproof

\begin{prop}The subalgebra $\A({\sf M}^u_{\theta;k})$ of coinvariants in $\B({\sf M}_{\theta;k})\lcross
H_F$ for the coaction $\delta_u$ is generated by elements of the
form $M^j_{ab}\otimes u_j$.
\end{prop}

\proof Clearly one has for all group-like elements $h\in H_F$ that
$$\delta_u(M^j_{ab}\otimes h)=u_j^*h\otimes M^j_{ab}\otimes h,$$
whence for coinvariants we need to take $h=u_j$.\endproof

We can explicitly compute the relations between generators of the
algebra $\A({\sf M}^u_{\theta;k})$ using the algebra relations of
$\B({\sf M}_{\theta;k})\lcross H_F$, obtaining
$$(M^j_{ab}\otimes u_j)(M^l_{cd}\otimes u_l)=(M^l_{cd}\otimes
u_l)(M^j_{ab}\otimes
u_j)F^{-2}(u_l,\tau_j^*)F^{-2}(\tau_l,\tau_j)F^{-2}(\tau_l^*,u_j) .$$
%
In particular, we can
take $(r_1,r_2)\in \ZZ^2$ to be a pair of integers and set
$$u=(u_1,u_2,u_3,u_4):=(\tau_1^{m_1},\tau_2^{m_2},\tau_3^{m_3},\tau_4^{m_4})$$
with $(m_1,m_2,m_3,m_4):=(r_1,r_1,r_2,r_2)$, as we did in the charge
one case. For such $u$, the commutation relations in $\A({\sf
M}^u_{\theta;k})$ reduce to
$$
(M^j_{ab}\otimes u_j)(M^l_{cd}\otimes
u_l)=\eta_{jl}^{m_l+m_j-1}(M^l_{cd}\otimes u_l)(M^j_{ab}\otimes
u_j).
$$
Let us check that there is a choice of integers
$r_1$, $r_2$ for which the parameter space $\A({\sf
M}^u_{\theta;k})$ is commutative. It is easy to see from
Eq.~\eqref{eqn eta matrix} that, whenever both $j,l\in\{1,2\}$ or
$j,l\in\{3,4\}$, the deformation parameter $\eta_{jl}$ is
automatically equal to $1$ and so these generators always commute.
Without loss of generality we consider the non-trivial case
$j\in\{1,2\}$ and $l\in\{3,4\}$, where the corresponding generators
fail to commute by a factor of $\eta_{jl}^{m_l+m_j-1}$. By
assumption we have that $m_l=r_1$ and $m_j=r_2$, so it follows that
any choice of $r_1$, $r_2$ for which $r_1+r_2=1$ makes the resulting
algebra $\A({\sf M}^u_{\theta;k})$ commutative.

Of course, for these parameter spaces there is a great deal of gauge
freedom left to be removed. As shown in \cite{bl:adhm}, the ADHM
construction does not depend on the choice of bases for the vector
spaces $\mH$, $\mL$ in the monad \eqref{eqn module monad}, whereas
making a unitary change of basis of $\mK$ which respects the
self-conjugacy property of the monad ({\em i.e.} acting with an
element of the unitary group $\Sp(\mK)$) results in a projection
which is unitarily equivalent to $\widetilde\Pp$. Removing the extra
gauge parameters corresponding to these degrees of freedom ought to
be straightforward, since the computation is entirely classical,
although we postpone an explicit computation to future work.

\subsection*{Acknowledgments}
This work was partially funded by the `Italian project Cofin06 -
Noncommutative Geometry, Quantum Groups and Applications'. We are
grateful to Lucio Cirio and Richard J. Szabo for extremely helpful
conversations.

\appendix

\section{Quantum Families of Maps}\label{families}In this appendix
we briefly review the notion of representability of functors and the
corresponding notion of universal objects. These are
of paramount importance in the present article, since they
tie together the various notions of universality that we use.

Let $\mathfrak{C}$ be a (locally small) category; for each pair of
objects $A$, $B$ of $\mathfrak{C}$, we write $\textup{Mor}(A,B)$ for
the set of morphisms from $A$ to $B$. Let
$\mathcal{F}:\mathfrak{C}\to \textsf{Set}$ be a covariant functor
from $\mathfrak{C}$ to the category $\textsf{Set}$ of sets.

\begin{defn}A {\em representation} of the functor $\mathcal{F}$ is a pair
$(M,\Phi)$, where $M$ is an object of $\mathfrak{C}$ and
$\Phi:\textup{Mor}(M,-)\to \mathcal{F}$ is an
isomorphism of functors ({\em i.e.} a natural transformation whose
component morphisms are all isomorphisms).  If such a representation
$(M,\Phi)$ exists, then the functor $\mathcal{F}$ is said to be {\em
representable}.
\end{defn}

From Yoneda's lemma one knows that natural transformations from
$\textup{Mor}(M,-)$ to $\mathcal{F}$ are in bijective correspondence
with elements of $\mathcal{F}(M)$ \cite{maclane}. Indeed, given a
natural transformation $\Phi:\textup{Mor}(M,-)\to \mathcal{F}$,
there is a corresponding element $\sigma\in\mathcal{F}(M)$ defined
by $\sigma:=\Phi_M(\id_M)$. Conversely, given $\sigma\in
\mathcal{F}(M)$, we can define a natural transformation
$\Phi:\textup{Mor}(M,-)\to \mathcal{F}$ by
$$\Phi_X(\delta):=(\mathcal{F}\circ\delta)(\sigma),\qquad \textup{for} \quad \delta\in
\textup{Mor}(M,X).$$
This leads to the following definition.
\begin{defn}A {\em universal object} for the functor $\mathcal{F}$ is a pair
$(M,\sigma)$, where $M$ is an object of $\mathfrak{C}$ and $\sigma$
is an element of the set $\mathcal{F}(M)$ with the property that for
every pair $(Y,\nu)$ with $Y$ an object of $\mathfrak{C}$ and $\nu$
an element of $\mathcal{F}(M)$, there is a unique morphism
$\Lambda\in\textup{Mor}(M,Y)$ such that $(\mathcal{F}\circ
\Lambda)(\sigma)=\nu$.\end{defn}

From the above argument it follows that representations of
$\mathcal{F}$ are in one-to-one correspondence with universal
objects  for $\mathcal{F}$. Of course, it is not necessarily the
case that a functor is representable, but if so, the
corresponding universal object is unique up to a unique isomorphism.
This abstract categorical set-up is extremely useful when applied to
the following examples.

\begin{example}First we recall the instanton moduli functor
$\mathcal{F}:\textsf{Alg}\to \textsf{Set}$ defined in
Remark~\ref{def mod}, which assigns to each unital $*$-algebra $A$
the set $\mathcal{F}(A)$ of equivalence classes of families of
instantons parameterised by $A$. To be more precise, we can define a
functor $\mathcal{F}_k$ by considering only families of instantons
with a fixed topological charge $k$. A (fine) moduli space of charge
$k$ instantons is a universal object representing the functor
$\mathcal{F}_k$.

Clearly, this set-up is usually far too naive for such a moduli
space to exist even in the classical case, but this example is
sufficient to illustrate why one should allow for the possibility of
noncommutative moduli spaces. The moduli space is necessarily an
object in the source of the functor $\mathcal{F}_k$ so that,
when allowing noncommutative parameter spaces,
one also needs to allow for the possibility of noncommutative moduli
spaces.\end{example}

\begin{example}Let $\mathfrak{C}$ be the
category whose objects are unital $C^*$-algebras. For any two
objects $A$ and $B$ the set of morphisms $\textup{Mor}(A,B)$ is the
set of all non-degenerate $*$-homomorphisms from $A$ to $B$.
Fix a pair of objects $A$, $B$ of $\mathfrak{C}$ and
define a functor $\mathcal{F}:\mathfrak{C}\to\textsf{Set}$ by
setting
$$\mathcal{F}(C):=\textup{Mor}(B,C\otimes A),$$ {\em i.e.} we assign to each $C^*$-algebra $C$
the set of all morphisms $\delta_C:B\to C\otimes A$. We say  \cite{sw:proc} that
$\delta_C$ is a {\em quantum family of maps} labeled by $C$.

In this case a universal object for $\mathcal{F}$ is a $C^*$-algebra
$M$ equipped with a morphism $\delta\in\textup{Mor}(B,M\otimes A)$
such that, for any $C^*$-algebra $C$ and any quantum family of maps
$\delta_C$ labeled by $C$, there exists a unique morphism
$\Lambda\in \textup{Mor}(M,C)$ and the
diagram
$$\begin{CD} B @>\delta>> M\otimes A
\\ @VV\id V @VV\Lambda\otimes\id V \\ B @>\delta_C>> C\otimes A
\end{CD}$$
is commutative. When $A,B,C$ are commutative
$C^*$-algebras, there exist compact Hausdorff topological spaces
$\Omega_A$, $\Omega_B$, $\Omega_C$ such that $A=C(\Omega_A)$  and so
on. A morphism $\delta_C\in \textup{Mor}(B,C\otimes A)$ corresponds
to a continuous map $\Omega_A\times\Omega_C\to \Omega_B$, {\em i.e.}
a continuous family of maps from $\Omega_A$ to $\Omega_B$
parameterised by the space $\Omega_C$. As explained in
\cite{sol:qfm}, the universal object corresponds to the space of
{\em all} continuous maps from $\Omega_A$ to $\Omega_B$. The
situation where the $C^*$-algebras are noncommutative
is a natural generalisation of this.\end{example}

As we have seen in the present article, one does not always need to
consider $C^*$-completions and can usually work perfectly well at
the level of $*$-algebras (this is also the usual setting for the
algebraic theory of quantum groups \cite{ma:book}). As usual,
we think of a noncommutative $*$-algebra $A$ as
the algebras of coordinate functions on some underlying `noncommutative
space' $\Omega_A$, with the $*$-structure interpreted as viewing
$\Omega_A$ as a `real form' of some complex affine space (see
\cite{bm:qtt} for further discussion in this direction). Of course,
one can add more structure if one wishes, such as requiring
$A$ to be Noetherian if one wants something resembling a
`$*$-algebraic variety', although we shall be deliberately vague
about this point.

\begin{example}Let $\mathfrak{C}$ be the category
of unital $*$-algebras, with morphisms given by non-degenerate
$*$-homomorphisms. The general principle of the previous example
still applies, now in the setting of $*$-algebraic geometry. Given a
pair $A$, $B$ of objects in the category, an element $\delta_C\in
\textup{Mor}(B,C\otimes A)$ is a quantum family of maps from
$\Omega_A$ to $\Omega_B$ parameterised by the noncommutative `space'
$\Omega_C$.\end{example}

\begin{example}In the situation of the previous example, we set $B=A$ and define a functor $\mathcal{F}:\mathfrak{C}\to
\textsf{Set}$ by assigning to each object $C$ the set of all
non-degenerate $*$-algebra maps $\delta_C:A\to C\otimes A$. One may
show \cite{sol:qfm} that the universal object $(M,\delta)$ is
automatically a bialgebra, whose coproduct and counit we denote
$\Delta_M$, $\ep_M$, and that it obeys the additional properties
$$(\id\otimes\delta)\circ\delta=(\Delta_M\otimes\id)\circ\delta,\qquad
(\ep_M\otimes \id)\circ\delta=\id,$$ {\em i.e.} $\delta$ makes $A$
into a left $M$-comodule algebra. We say that a pair $(C,\delta_C)$
obeying these properties is a {\em transformation bialgebra} for the
algebra $A$. The universal object is called the {\em universal
transformation bialgebra} \cite{lprs:ncfi}. In the classical case,
it just corresponds to the semigroup of all algebraic maps from the
commutative space $\Omega_A$ to itself.\end{example}

\begin{example}Let $H$ be a coquasitriangular Hopf $*$-algebra and take
$\mathfrak{C}$ to be the category of left $H$-comodule algebras
which, as discussed in \S\ref{section hopf algebra prelims}, is a
braided monoidal category. Once again we fix an algebra $A$ in the
category, but we take now the functor
$\mathcal{F}:\mathfrak{C}\to \textsf{Set}$ to be the one which
assigns to each object $C$ of $\mathfrak{C}$ the set of all {\em
braided} morphisms $\delta_C:A\to C\utimes A$, where $\utimes$ is
the tensor product induced by the braiding. The universal object is
the {\em universal braided transformation bialgebra} for $A$.

This is the strategy we adopted in \S\ref{section quantum conformal
group}, where we took $H=\A(\TT^2)$ and $A=\A(\C^4)$. With the
additional requirement that the bialgebra must respect the
quaternionic structure of $\HH^2\simeq\C^4$ (hence inducing the
$*$-structure in Eq.~\eqref{eqn defining M(H)} as in
\cite{lprs:ncfi}), we found that the universal transformation
bialgebra in the category is the matrix bialgebra $\A(\M(2,\HH))$.
The fact that the quantisation functor is an isomorphism of braided
monoidal categories means that it preserves the universality
property, so that now viewing $\mathcal{F}$ as a functor from the
category of left $H_F$-comodule algebras to the category of sets,
the universal transformation bialgebra for $\A(\C^4_\theta)$ is the
braided matrix bialgebra $\B(\M_\theta(2,\HH))$ of  \S\ref{section
braided groups}.\end{example}

Our final example concerns the construction of parameter spaces of
module maps as universal objects. It is more general than the
previous examples, which considered algebra maps, but it still uses
a universality property to define the `space of all maps'. The
example illustrates that if one changes the source category of a
functor then the problem of representability can alter dramatically.
We stress once again that in looking for moduli spaces of
instantons, our philosophy is to look not for a set of objects but
rather for a space which parameterises those objects, that is to say
we ask for some geometric structure. In categorical terms, this
means defining a functor from the category of algebras to the
category of  sets and then looking for the moduli space as a
universal object, which is by definition an object in the source
category and so necessarily an algebra.

\begin{example}
In \S\ref{nc monads}, we
considered right module maps
$$\sigma_z:\mH\otimes\A(\C^4)(-1)\to\mK\otimes\A(\C^4)$$ which are linear in the generators $z_1,\ldots,z_4$ of
$\A(\C^4)$ and then looked for the space of all such maps. To view
this in a categorical setting, we take $\mathfrak{C}$ {\em a priori}
to be the category of unital algebras ($*$-structures
are not required at this stage) and consider the functor
$\mathcal{F}:\mathfrak{C}\to \textsf{Set}$ which assigns to each
algebra $C$ the set of all right $\A(\C^4)$-module maps
$$\delta_C:\mH\otimes\A(\C^4)(-1)\to C\otimes\mK\otimes\A(\C^4)$$
which are linear in the generators $z_1,\ldots,z_4$ of $\A(\C^4)$.
We would like to find the space of all maps $\sigma_z$ in terms of a
universal object for this functor ({\em i.e.} by proving that it is
representable).
Following the approach taken in \S\ref{nc monads} and in the above
examples, we try and prove representability in this case by
explicitly constructing the universal algebra. It is straightforward
to see that the universal algebra, if it exists, must be generated
by the elements $M^{\alpha}_{ab}$ which define the map
\begin{equation}\sigma_z:u_b\otimes
Z\mapsto\sum_{a,\alpha} M^\alpha_{ab}\otimes v_a\otimes z_\alpha
Z,\qquad Z\in \A(\C^4) ,
\end{equation} where $j=1,\ldots,4$ and
$a=1,\ldots,2k+2$, $b=1,\ldots,k$. In our approach, we need to find
an algebra structure on this set of functions $M^j_{ab}$. However,
the construction fails at this point: the objects
$\mH\otimes\A(\C^4)$ and $\mK\otimes\A(\C^4)$ are only
$\A(\C^4)$-modules and do not themselves have an algebra structure, so there is nothing
to determine an algebra structure on the matrix elements $M^j_{ab}$
and one has to make a choice. 

One could alternatively consider looking for the set of all module maps simply as a vector space, rather than looking for its coordinate algebra. However, this is not very natural as it does not imply any geometric structure; also there
does not seem to exist a corresponding notion of universality.

One way to proceed is to look for the space of all such maps
$\sigma_z$ as a classical object, just as we did in \S\ref{nc
monads}. This means taking the source category $\mathfrak{C}$ of the
functor $\mathcal{F}$ to be the category of {\em commutative} unital
algebras: it is perfectly natural to assume in this way that the
algebra $\A(\M(\mH,\mK))$ generated by the coordinate functions
$M^j_{ab}$ is commutative, hence giving the space $\M(\mH,\mK)$ the
structure of a classical algebraic variety. By restricting the
functor in this way, it  becomes representable
with $\A(\M(\mH,\mK))$ as the universal object.

Now that we have constructed $\A(\M(\mH,\mK))$ as a suitable
parameter space of maps, we proceed just as in \S\ref{nc monads} to
show that, in fact, this algebra is an object in the category of
left $H$-comodules. Viewed in this way, the noncommutative parameter
spaces that we construct are just the canonical deformations of the
corresponding classical objects and so, in this sense, they are the
most natural objects to work with.
\end{example}

\end{document}